\documentclass[12pt]{article}
\usepackage{a4}
\usepackage{amsmath}
\usepackage{amssymb}
\usepackage{amsbsy}
\usepackage[ctr2,eng]{myarticl}
\usepackage[frame,curve,knot]{xy}
\usepackage{ifthen}
\usepackage{epsfig}

\newcommand{\mc}{\mathcal}
\newcommand{\A}{\mathcal{A}}
\newcommand{\Anull}{\A ^\circ }
\newcommand{\Gam}{\varGamma }
\newcommand{\x}{X}
\newcommand{\w}{\theta }
\newcommand{\cconj}[1]{#1^\mathrm{c}}
\newcommand{\cb}{\eta }
\newcommand{\crs}{\Gamma }
\newcommand{\tors}{T}
\newcommand{\szushg}{\nabla _{\!\!\mc{S}}}
\newcommand{\multw}{{\mathrm{m}_\wedge }}
\newcommand{\multcl}{\mathrm{m}_{\scriptscriptstyle \mathrm{Cl}}}
\newcommand{\multcls}{\mathrm{m}_{\scriptscriptstyle \mathrm{Cl,S}}}
\newcommand{\spinsig}[1]{\nu _{#1}}
\newcommand{\invcsu}{\dagger }
\newcommand{\invncsu}{\star }
\newcommand{\invslR}{\sharp }

\newcommand{\qp}{[2]}
\newcommand{\qm}{\hat q}
\newcommand{\oC}{\check C}
\newcommand{\uC}{\hat C}
\newcommand{\ts}{\tilde{\sigma }}
\newcommand{\real}{\mathbb{R}}
\newcommand{\realx}{\mathbb{R}^\times }
\newcommand{\comp}{\mathbb{C}}
\newcommand{\compx}{\mathbb{C}^\times }
\newcommand{\ot}{\otimes }
\newcommand{\otA}{{\otimes _{\scriptscriptstyle \A }}}
\newcommand{\SLq}[1]{\mathrm{SL}_q(#1)}
\newcommand{\SUq}[1]{\mathrm{SU}_q(#1)}
\newcommand{\OSLq}[1]{\mc{O}(\SLq{#1})}
\newcommand{\OSUq}[1]{\mc{O}(\SUq{#1})}
\newcommand{\Uqsl}[1]{U_q(\mathrm{sl}_{#1})}
\newcommand{\tUqsl}[1]{\tilde{U}_q(\mathrm{sl}_{#1})}
\newcommand{\Gamw}[1][]{\Gam ^{\wedge #1}}
\newcommand{\Gamt}[1][]{\Gam ^{\otA #1}}
\newcommand{\ctr}[2]{\bigl\langle \!\!\langle #1,#2\rangle \!\!\bigr\rangle }
\newcommand{\hmet}[2]{\langle #1,#2\rangle }
\newcommand{\Cl}{\mathrm{Cl}}
\newcommand{\lact}{\triangleright }
\newcommand{\ract}{\triangleleft }
\newcommand{\koun}{\varepsilon }
\newcommand{\kopr}{\varDelta }
\newcommand{\lkow}{\varDelta _{\scriptscriptstyle \mathrm{L}}}
\newcommand{\rkow}{\varDelta _{\scriptscriptstyle \mathrm{R}}}
\renewcommand{\linv}[1]{#1_{\scriptscriptstyle \mathrm{L}}}
\renewcommand{\rinv}[1]{#1_{\scriptscriptstyle \mathrm{R}}}
\newcommand{\tr}{\mathrm{Tr}}
\newcommand{\lfunc}{{\boldsymbol \ell}}
\newcommand{\tlfunc}{\tilde{\boldsymbol \ell}}
\newcommand{\Xalg}[1][]{{<}\mc{X}_{#1}{>}}
\newcommand{\coalg}{{\boldsymbol u}}
\newcommand{\Cas}{\mc{Z}_q}
\newcommand{\mybox }{\rule[-.5ex]{0pt}{2ex}}

\title{Spin Geometry on Quantum Groups via Covariant Differential
Calculi\thanks{This paper was supported by the Deutsche
Forschungsgemeinschaft}}
\author{Istv\'an Heckenberger\thanks{e-mail:
heckenbe@mathematik.uni-leipzig.de}\\
{\small Universit\"at Leipzig, Augustusplatz 10-11,}\\
{\small 04109 Leipzig, Germany}}
\date{}

\begin{document}
\sloppy

\maketitle

\begin{abstract}
Let $\A $ be a cosemisimple Hopf $*$-algebra with antipode $S$ and let $\Gam $
be a left-covariant first order differential $*$-calculus over $\A $ such that
$\Gam $ is self-dual (see Section \ref{sec-diffcalc}) and invariant under the
Hopf algebra automorphism $S^2$.
A quantum Clifford algebra $\Cl (\Gam ,\sigma ,g)$ is introduced which
acts on Woronowicz' external algebra $\Gamw $. A minimal left ideal of
$\Cl (\Gam ,\sigma ,g)$ which is an $\A $-bimodule is called a spinor module.
Metrics on spinor modules are investigated. The usual notion of a linear left
connection on $\Gam $ is
extended to quantum Clifford algebras and also to spinor modules. The
corresponding Dirac operator and connection Laplacian are defined.
For the quantum group $\SLq 2$ and its bicovariant $4D_\pm $-calculi
these concepts are studied in detail. A generalization of Bochner's theorem
is given. All invariant differential operators over a given spinor module are
determined. The eigenvalues of the Dirac operator are computed.
\end{abstract}

\noindent
{\small
}

\setcounter{section}{-1}
\section{Introduction}

Noncommutative geometry has been invented in the eighties as a new field
by the pioneering work of A.~Connes \cite{b-Connes1}.
Nowadays it is commonly expected that it will provide new mathematical tools
and methods for applications in theoretical physics
(quantum gravity, physics at small distances).
On the other hand, quantum groups also arose in the eighties as new classes
of noncommutative non-cocommutative
Hopf algebras \cite{a-FadResTak1}.
A general framework for bicovariant noncommutative differential
calculi over Hopf algebras has been provided by S.\,L.~Woronowicz
\cite{a-Woro2}. In the meantime a deep algebraic theory of such calculi
has been developped (see, for instance, \cite{b-KS}, Chapter 14).
Unfortunately, the bicovariant differential calculus on quantum groups does
not fit into the realm of the noncommutative geometry of A.~Connes.
As a first step to connect both theories one needs spin structures
and Dirac operators on quantum groups. In this paper we study whether
covariant differential calculi permit the setup of a spin geometry.
The bridge to A.~Connes' noncommutative geometry will be considered in a
forthcoming paper.

The theory of spin and related structures for (pseudo)riemannian manifolds
has its origin in Dirac's relativistic theory of the electron (1928).
Comprehensive introductions into the subject can be found (for instance)
in the books by E.~Cartan \cite{b-Cartan66}, A.~Crumeyrolle \cite{b-Crumey90},
and H.\,B.~Lawson and M.-L.~Michelsohn \cite{b-LawMich89}.
Recently some attempts have been made to generalize these structures to
noncommutative algebras. Connections on bimodules of differential one-forms
have been examined by several authors \cite{a-CQ1}, \cite{a-DMMM95},
\cite{a-BDMS1}, \cite{a-DHLS}, \cite{a-HeckSchm1}. In the paper of
R.~Bautista et.\,al.\ \cite{a-BCDRV96} quantum Clifford algebras are
constructed in the case where the braiding on the tensor product of one-forms
satisfies the Hecke relation.

The aim of this paper is to present generalizations of certain
well-known structures of the spin geometry to left-covariant differential
$*$-calculi on quantum groups.

Throughout we deal with a cosemisimple Hopf $*$-algebra $\A $ and with a
self-dual $S^2$-invariant left-covariant differential $*$-calculus over $\A $.
A definition of the corresponding quantum Clifford algebra and its main
involution is given. It is proved that the quantum Clifford algebra has a
representation on Woronowicz' external algebra.
Spinor modules are defined as left ideals and $\A $-subbimodules of the
quantum Clifford algebra. Conditions for the extension
of a connection on the differential calculus to the latter objects are given.
The connection Laplacian and the Dirac operators are constructed.
The symmetry of the Dirac operator with respect to metrics on
spinor modules is analyzed. 
After a discussion of the general situation the quantum group $\SLq 2$
(with three non-isomorphic $*$-structures) and the two 4-dimensional
bicovariant differential calculi on it are studied in detail.
It is proved that there exists a minimal spinor module $\mc{S}_0$ of the
quantum Clifford algebra. We introduce a Hopf $*$-algebra
$\tilde{\A }$ which contains $\A =\OSLq 2$ as a Hopf subalgebra.
Then one can define a right coaction of $\tilde{\A }$ on $\mc{S}_0$ which
is compatible with the
Clifford multiplication and the original right coaction of the differential
calculus. Therefore $\tilde{\A }$
can be interpreted as the function algebra of a covering of the quantum group
$\SLq 2$. Also the uniqueness of the metric on $\mc{S}_0$ is proved.
It is important to emphasize that this metric is hermitean and non-degenerate
but it is not positive definite.
Connections on the quantum Clifford algebra and on $\mc{S}_0$ are found.
It turns out that with our definitions there exists no torsion free connection
on the differential one-forms. Finally all invariant differential operators
on the spinor module $\mc{S}_0$ are determined. The eigenvalues of the
Dirac operator are computed and a generalization of Bochner's
theorem is proved.

The paper is organized in the following way. In Section \ref{sec-coshsa} the
definition and some properties of cosemisimple Hopf $*$-algebras are given.
In Section \ref{sec-diffcalc} we collect the assumptions on the
differential calculus over the Hopf $*$-algebra. In Section \ref{sec-cliff} we
introduce and study quantum Clifford algebras, spinor modules
and metrics on it. In Section \ref{sec-connection} we investigate
linear connections on differential forms, quantum exterior
and quantum Clifford algebras and on spinor modules. We introduce a
duality of left connections which makes the definition of the connection
Laplacian possible. In Section \ref{sec-Dirac} the Dirac operator is defined
and some of its properties are proved. In Section \ref{sec-example}
we apply the general theory developped in the preceding sections to the
4-dimensional
bicovariant differential calculi on the quantum groups
$\SUq 2$, $\SUq{1,1}$ and $\SLq{2,\real }$. All structures are worked out
in great detail. Quantum spin groups for the differential
structures on these quantum groups and metrics and connections
on the spinor module $\mc{S}_0$ are determined.
In Section \ref{sec-invdo} we use generalized $\ell $-functionals and
Clebsch-Gordan coefficients to obtain all invariant differential operators
on the spinor module $\mc{S}_0$ and on one of its subbimodules $\mc{S}_0^+$.
One of our main results (Theorem \ref{t-invdosp}) states that the algebra of 
invariant differential operators on $\mc{S}_0^+$ is commutative and generated
by two elements and one relation (\ref{eq-invdorel}). In Section
\ref{sec-invdo} we also compute the
eigenvalues of the Dirac operator and prove a generalization of Bochner's
theorem (Theorem \ref{t-SLq2Bochner}). In two appendices
some properties of the Hopf algebras $\OSLq 2$ and $\Uqsl 2$ and a
corresponding graphical calculus for the morphism spaces are recalled.

Let us fix some notation. If not otherwise stated, we follow the definitions
and conventions of the monograph \cite{b-KS}. In this paper $\A $ denotes
a cosemisimple Hopf $*$-algebra over the complex numbers with coproduct
$\kopr $, counit $\koun $ and antipode $S$. We use Sweedler's notation
$\kopr (a)=a_{(1)}\ot a_{(2)}$. The symbol $\Mor (v,w)$ denotes the vector
space of intertwiners $T$ of corepresentations $v$ and $w$ of $\A $.
Throughout we use Einstein's convention to sum over repeated indices.

The author is grateful to K.~Schm\"udgen for suggesting the subject of this
paper and H.\,B.~Rademacher for hints to clarify the commutative case.

\section{Cosemisimple Hopf $*$-algebras}
\label{sec-coshsa}

A Hopf algebra $\A $ is called a Hopf $*$-algebra if $\A $ is a $*$-algebra
with involution
$*:\A \to \A $ such that $\kopr (a^*)=(*\ot *)\kopr (a)$ for any $a\in \A $.
A Hopf algebra $\A $ is called cosemisimple, if any corepresentation of $\A $
is a direct sum of simple corepresentations. If not otherwise
stated, $\A $ always denotes a cosemisimple Hopf $*$-algebra.

One of the most important properties of cosemisimple Hopf algebras is the
existence of a unique
left- and right-invariant functional (the Haar functional) $h$ on $\A $
such that $h(1)=1$. Moreover, the Haar functional is left- and right-regular,
i.\,e.\ $b\in \A $ and $h(ab)=0$ ($h(ba)=0$) for any $a\in \A $ implies
$b=0$.

Since the functional $h'$ on $\A $, defined by $h'(a):=\overline{h(a^*)}$
for all $a\in \A $ is also left- and right-invariant and $h'(1)=1$, we have
$h'=h$ and therefore $h(a^*)=\overline{h(a)}$ for any $a\in \A $.

It is well known that there is an automorphism $\rho $ of $\A $ such that
$h(ab)=h(b\rho (a))$ for any $a,b\in \A $. More precisely, if $w=(w_{ij})$ is
a simple corepresentation of $\A $ and
$F=(F_{kl})\in \Mor (w,w^{\contr \contr })$, $F\not=0$,
then $F$ is an invertible matrix with complex entries, $\tr F\not= 0$,
$\tr F^{-1}\not= 0$ and
\begin{align}\label{eq-rhow}
\rho (w_{kl})&=\frac{F_{km}F_{nl}\tr F^{-1}}{\tr F}w_{mn}.
\end{align}
Let us define a mapping $\beta :\A \to \A $ by $\beta (a):=\rho (a)^*$.
Then we have $h(ab)=h(b\beta (a)^*)$ for any $a,b\in \A $ and
$\beta $ becomes an algebra involution of $\A $. Indeed, $\beta $ is an
antilinear antiendomorphism of $\A $ and
\begin{align}
h(ab)=h(b\beta (a)^*)=\overline{h(\beta (a)b^*)}
=\overline{h(b^*\beta ^2(a)^*)}=h(\beta ^2(a)b)
\end{align}
for any $a,b\in \A $. Hence $a=\beta ^2(a)$ for any $a\in \A $.

\begin{lemma}
For the algebra involution $\beta $ of $\A $ the following formula holds:
\begin{align}\label{eq-koprbeta}
\kopr (\beta (a))=\beta (a_{(1)})\ot S^2(a_{(2)}^*)\qquad \text{for all
$a\in \A $.}
\end{align}
\end{lemma}

\begin{bew}
Both sides of (\ref{eq-koprbeta}) are antilinear mappings from
$\A $ to $\A \ot \A $. Hence it suffices to prove the lemma for
$a=w_{ij}$, where $w_{kl}$ are matrix elements of an irreducible
corepresentation of $\A $ such that $\kopr (w_{ij})=w_{ik}\ot w_{kj}$.
Let $F$ be an invertible morphism of the corepresentations $w$ and
$w^{\contr \contr }$. For the left hand side of (\ref{eq-koprbeta}) we obtain
\begin{align*}
\kopr (\beta (w_{ij}))&=\kopr (\rho (w_{ij})^*)
=\kopr \left(\frac{\overline{F_{ik}}\,\overline{F_{lj}}\,
\overline{\tr F^{-1}}}{\overline{\tr F}}w^*_{kl}\right)
\end{align*}
by (\ref{eq-rhow}).
Since $\kopr $ is a $*$-homomorphism, the latter is equal to
\begin{align*}
\frac{\overline{F_{ik}}\,\overline{F_{lj}}\,
\overline{\tr F^{-1}}}{\overline{\tr F}}w^*_{km}\ot w^*_{ml}
&=\frac{\overline{F_{ik}}\,\overline{F_{mn}}\,
\overline{\tr F^{-1}}}{\overline{\tr F}}w^*_{km}\ot 
\overline{F^{-1}_{nr}}\,\overline{F_{lj}}w^*_{rl}\\
&=\beta (w_{in})\ot (F^{-1}_{nr}w_{rl}F_{lj})^*.
\tag{$*$}
\end{align*}
Finally, $F\in \Mor (w,w^{\contr \contr })$ gives
$w_{rl}F_{lj}=S^{-2}(S^2(w_{rl})F_{lj})=S^{-2}(F_{rl}w_{lj})$.
Hence $(F^{-1}_{nr}w_{rl}F_{lj})^*=(F^{-1}_{nr}F_{rl}S^{-2}(w_{lj}))^*
=S^2(w_{nj}^*)$ and (\ref{eq-koprbeta}) follows from ($*$).
\end{bew}

\section{$S^2$-invariant differential calculi}
\label{sec-diffcalc}

Let $\A $ be a cosemisimple Hopf $*$-algebra and $(\Gam ,\dif )$ a
left-covariant first order differential $*$-calculus over $\A $.
Let $\mc{X}$ denote the quantum tangent space of $\Gam $.
Since $S^2$ is a Hopf algebra automorphism of $\A $, there is a left-covariant
first order differential calculus $(\Gam ',\dif ')$ with quantum tangent space
$\mc{X'}:=S^2(\mc{X})$ (see \cite{a-Heck00a}). We call $(\Gam ,\dif )$
$S^2$-\textit{invariant}, if $(\Gam ,\dif )$ and $(\Gam ',\dif ')$ are
isomorphic, i.\,e.\ $\mc{X'}=\mc{X}$.
\textbf{
In this article all first order differential calculi are assumed to be 
$S^2$-invariant left-covariant differential $*$-calculi.
}

Let $\{\x _i\,|\,i=1,\ldots ,n=\dim \mc{X}\}$ be a basis of $\mc{X}$ and
let $\{\w _i\,|\,i=1,\ldots ,n\}$ be the dual basis of the vector space
of left-invariant 1-forms $\linv{\Gam }$.
This means that for the differential mapping
$\dif :\A \to \Gam $ the formula
\begin{align}\label{eq-diffmap}
\dif a&=a_{(1)}\x _i(a_{(2)})\w _i,\qquad a\in \A ,
\end{align}
holds.
Recall that $\x _i\in \Anull $ and there exist functionals
$f^i_j\in \Anull $ such that
$\kopr \x _i=\koun \ot \x _i +\x _j\ot f^j_i$. Moreover, the formula
$\w _ia=a_{(1)}f^i_j(a_{(2)})\w _j$ holds for $a\in \A $ and $i=1,\ldots ,n$.
If $S^2(\x _i)=F^i_j\x _j$ for some $F^i_j\in \mathbb{C}$, then the mapping
$\w _i\mapsto S^2(\w _i):=F^j_i\w _j$ extends uniquely to an isomorphism
of the left-covariant $\A $-bimodule $\Gam $ and we have
$S^2(\dif a)=\dif S^2(a)$ for all $a\in \A $.

Recall that the involution $*$ of $\Anull $ is defined by
\begin{gather}\label{eq-stern}
f^*(a)=\overline{f\bigl(S(a)^*\bigr)}\quad \text{or}\quad
f(a^*)=\overline{S(f)^*(a)}.
\end{gather}
Since $(\Gam ,\dif )$ is a
$*$-calculus, we have $\mc{X}^*=\mc{X}$.
Hence there is a matrix $E=(E^i_j)$ such that $\x ^*_i=E^i_j\x _j$
(and therefore $\w ^*_i=-\w _j\overline{E^j_i}$) for
$i=1,\ldots ,n$. Let now $B=(B^i_j)_{i,j=1,\ldots ,n}$
denote the complex matrix $B=\bar{F}E$. This implies that
\begin{align}
S^2(\x _i)^*=B^i_j\x _j.
\end{align}

\begin{satz}\label{s-betaGam}
The setting $\beta (\w _i):=\overline{B^j_i}\w _j$, $\beta (a)$ as
in the previous
section, defines an involution of the $\A $-bimodule $\Gam $. Moreover,
$\beta (\dif a)=-\dif \beta (a)$ for any $a\in \A $.
\end{satz}

\begin{bew}
Since $\x _i$ is an element of the Hopf algebra $\Anull $, from
$S(S(\x _i)^*)^*=\x _i$ we obtain
$S^2(S^2(\x _i)^*)^*=\x _i$ for any $i=1,\ldots ,n$. Therefore
\begin{align}
\x _i=S^2(B^i_j\x _j)^*=\overline{B^i_j}S^2(\x _j)^*
=\overline{B^i_j}B^j_k\x _k
\end{align}
and so $\overline{B}B=\id $. Hence
$\beta ^2(\w _i)=\beta (\overline{B^j_i}\w _j)
=B^j_i\overline{B^k_j}\w _k=\w _i$ for any $i=1,\ldots ,n$.
Now we only have to show that $\beta $ is well-defined,
that is $0=\beta (a)\beta (\w _i)-\beta (\w _j)\beta (a_{(1)}f^i_j(a_{(2)}))$
(which is formally the image under $\beta $ of the element
$\w _ia-a_{(1)}f^i_j(a_{(2)})\w _j$) for $a\in \A $ and $i=1,\ldots ,n$.
For this we compute
\begin{align}
\kopr (S^2(\x _i)^*)=&((S^2\ot S^2)\kopr (\x _i))^*
=1\ot S^2(\x _i)^*+S^2(\x _j)^*\ot S^2(f^j_i)^*\notag \\
=& 1\ot B^i_j\x _j+B^j_k\x _k\ot S^2(f^j_i)^*,\\
\kopr (B^i_j\x _j)=&B^i_j1\ot \x _j+B^i_j\x _k\ot f^k_j.
\end{align}
Hence we obtain $B^i_jf^k_j=B^j_kS^2(f^j_i)^*$,
i.\,e.\ $f^k_l=\overline{B^l_i}B^i_jf^k_j
=\overline{B^l_i}B^j_kS^2(f^j_i)^*$.
Therefore
\begin{align*}
&\beta (\w _j)\beta (a_{(1)}f^i_j(a_{(2)}))
=\overline{B^k_j}\w _k\beta (a_{(1)})
\overline{\overline{B^j_m}B^n_iS^2(f^n_m)^*(a_{(2)})}\\
&\quad =\overline{B^k_j}\beta(a_{(1)})f^k_l(S^2(a^*_{(2)}))\w _l
B^j_m\overline{B^n_i}S^2(f^n_m)(S^{-1}(a^*_{(3)}))\\
&\quad =\overline{B^n_i}\beta (a_{(1)})(S^2(f^k_l)S(f^n_k))(a_{(2)}^*)\w _l
=\overline{B^n_i}\beta (a_{(1)})\delta ^n_l\overline{\koun (a_{(2)})}\w _l
=\beta (a)\overline{B^l_i}\w _l.
\end{align*}
Here from the $1^\mathrm{st}$ to the $2^\mathrm{nd}$ line we used
(\ref{eq-koprbeta}) and (\ref{eq-stern}).
Similarly one can compute $\beta (\dif a)$.
\begin{align*}
&\beta (a_{(1)}\x _i(a_{(2)})\w _i)
=\beta (\w _i)\beta (a_{(1)})\overline{\x _i(a_{(2)})}
=\overline{B^j_i}\w _j\beta (a_{(1)})\overline{\x _i(a_{(2)})}\\
&\quad =\overline{B^j_i}\beta (a_{(1)})f^j_k(S^2(a^*_{(2)}))\w _k
S^{-2}(\x ^*_i)(S(a^*_{(3)}))\\
&\quad
=\overline{B^j_i}\beta (a_{(1)})(B^i_l\x _lS(f^j_k))(S(a^*_{(2)}))\w _k
=\beta (a_{(1)})S(\x _jS(f^j_k))(a^*_{(2)})\w _k\\
&\quad =\beta (a_{(1)})(-S^2(\x _k))(a^*_{(2)})\w _k
=-\beta (a_{(1)})\x _k(S^2(a^*_{(2)}))\w _k
=-\dif \beta (a).
\end{align*}
\end{bew}

If $(\Gam ,\dif )$ is a bicovariant FODC over $\A $ then there is a canonical
method (given by S.\,L.~Woronowicz \cite{a-Woro2})
to construct a differential calculus $(\Gamw ,\dif )$ over $\A $ with first
order part $(\Gam ,\dif)$. This construction is based on the existence
of a braiding, an automorphism $\sigma $ of the bicovariant
$\A $-bimodule $\Gam \otA \Gam $ satisfying the braid relation.
The mapping $\sigma $ is defined by the formula
\begin{align}\label{eq-bikovbraid}
\sigma (\rho \otA \rho ')&=\rho '\otA \rho,\qquad
\rho \in \linv{\Gam },\rho '\in \rinv{\Gam }.
\end{align}
However, for the existence of a braiding bicovariance is not necessary
(for an example see \cite{a-Heck00a}). Let $\multw $ denote the canonical
mapping
$\multw :\Gamt[2]\to \Gamw[2]=\Gamt[2]/\ker (\id -\sigma )$.

Let $(\Gam ,\dif )$ be a ($S^2$-invariant left-covariant) first order
differential $*$-cal\-cu\-lus over $\A $. Assume that there is an
automorphism $\sigma $ of the left-covariant $\A $-bimodule
$\Gam \otA \Gam $ such that
\begin{enumerate}
\item
$\sigma $ satisfies the braid relation
$\sigma _{12}\sigma _{23}\sigma _{12}=\sigma _{23}\sigma _{12}\sigma _{23}$
on $\Gamt[3]$,
\item \label{breq2}
$\ker (\id -\sigma )$ contains all elements
$\omega (a_{(1)})\otA \omega (a_{(2)})\in \Gam \otA \Gam $, where
$a\in \mc{R}_\Gam =\{b\in \A \,|\,\omega (b)=0,\koun (b)=0\}$,
\item
$S^2(\sigma \rho _2)=\sigma (S^2(\rho _2))$ for any $\rho _2\in \Gamt[2]$,\\
(where $S^2(\rho \otA \rho ')=S^2(\rho )\otA S^2(\rho ')$,
$\rho ,\rho '\in \Gam $,)
\item
$(\sigma \rho _2)^*=\sigma (\rho _2^*)$ or
$(\sigma \rho _2)^*=\sigma ^{-1}(\rho _2^*)$
for any $\rho _2\in \Gamt[2]$.
\end{enumerate}
Such a mapping $\sigma $ is called a \textit{braiding of} $\Gam $.
Now for any $k\in \mathbb{N}_0$ there is
an antisymmetrizer $A_k$, defined by $\sigma $ \cite{a-Woro2},
which is an automorphisms of the left-covariant
$\A $-bimodule $\Gamt[k]=\Gam \otA \cdots \otA \Gam $ ($k$ times).
In particular, we have $A_1=\id $ and $A_2=\id -\sigma $.
Moreover, the direct sum of the
kernels of the antisymmetrizer $A_k$ is a twosided ideal in
$\bigoplus _{k=0}^\infty \Gamt[k]$.
If we now define $\Gamw[k]:=\Gamt[k]/\ker A_k$,
$\Gamw:=\bigoplus _{k=0}^\infty \Gamw[k]$, then there is a canonical
differential calculus $(\Gamw ,\dif )$ over $\A $ with first order part
$(\Gam ,\dif )$.

Let us call a FODC $(\Gam ,\dif )$ \textit{self-dual}, if
there exists a braiding $\sigma $ of $\Gam $ and
a left-covariant $\sigma $-metric $g:\Gam \otA \Gam \to \A $ \cite{a-Heck99}.
More exactly, the mapping
$g$ should be a homomorphism of left-covariant $\A $-bimodules,
be non-degenerate, and should satisfy the equations
$g\sigma =g$ on $\Gamt[2]$
and $g_{12}\sigma ^\pm _{23}\sigma ^\pm _{12}=g_{23}$ on $\Gamt[3]$.
Further, for compatibility with the involution we assume that
\begin{align}\label{eq-gstar}
g(\rho \otA \rho ')^*&=g(\rho '{}^*\otA \rho ^*)\quad
\text{for $\rho ,\rho '\in \Gam $.}
\end{align}

\begin{lemma}
Let $\Gam $ be a bicovariant $\A $-bimodule with canonical braiding $\sigma $.
Suppose that there exists a homomorphism $g:\Gam \otA \Gam \to \A $
of bicovariant $\A $-bimodules. Then the equation $g\sigma =g$ is fulfilled
if and only if
\begin{align}\label{eq-gSquad}
g(S^2(\rho ')\otA \rho )&=g(\rho \otA \rho '),\qquad
\rho ,\rho '\in \linv{\Gam }.
\end{align}
\end{lemma}

\begin{bew}
First recall that
$S^2(\rho ')=S^2(\rho '_{(2)})\rho '_{(0)}S(\rho '_{(1)})$
for $\rho '\in \linv{\Gam }$, where
$\rkow ^2(\rho ')=\rho '_{(0)}\ot \rho '_{(1)}\ot \rho '_{(2)}$. Observe that
$\rho '_{(0)}S(\rho '_{(1)})\in \rinv{\Gam }$. Using
(\ref{eq-bikovbraid}), $g\sigma =g$ and left-covariance of $g$ it follows that
\begin{align*}
g(S^2(\rho ')\otA \rho )&=
S^2(\rho '_{(2)})g(\rho '_{(0)}S(\rho '_{(1)})\otA \rho )
=S^2(\rho '_{(2)})g\sigma (\rho \otA \rho '_{(0)}S(\rho '_{(1)}))\\
&=S^2(\rho '_{(2)})g(\rho \otA \rho '_{(0)})S(\rho '_{(1)})\\
&=S^2(\rho '_{(2)})S(\rho '_{(1)})g(\rho \otA \rho '_{(0)})
=g(\rho \otA \rho ')
\end{align*}
for $\rho ,\rho '\in \linv{\Gam }$.
The other direction of the assertion can be shown similarly.
\end{bew}

Now if $\Gam $ is an $S^2$-invariant left-covariant FODC then
we additionally require that left-covariant $\sigma $-metrics $g$
on $\Gam $ have to satisfy (\ref{eq-gSquad}).
This in turn implies that $g(S^2(\rho ))=S^2(g(\rho ))$ holds for
$\rho \in \Gam \otA \Gam $.
\textbf{
In this article we always assume that the first order differential calculus
$(\Gam ,\dif )$ is self-dual.
}

\begin{bem}
The property (\ref{eq-gSquad}) of the $\sigma $-metric $g$ will be used
only in the proof of Proposition \ref{s-laplso}.
However, the requirement $gS^2=S^2g$
is essential for the existence of the algebra involution $\beta $ of the
quantum Clifford algebra $\Cl (\Gam ,\sigma ,g)$, see Proposition
\ref{s-betaCl}.
\end{bem}

Suppose for a moment that $\A $ is a coquasitriangular Hopf algebra and
there is a simple corepresentation $u$ of $\A $ such that the
contragredient corepresentation $u\cont $ is isomorphic to $u$.
For example, this is the case for $\A =\OSLq 2$ and $u$ the fundamental
corepresentation of $\A $.
In \cite{a-Heck99}
it was proved that in this setting there exists a bicovariant
first order differential calculus $(\Gam (u),\dif )$ over $\A $ which is
self-dual. Moreover, there exists a bicovariant $\sigma $-metric on
$\Gam (u)$. This example will be considered in detail in Section
\ref{sec-example}.

\section{Quantum Clifford algebras}
\label{sec-cliff}

Let $(\Gam ,\dif )$ be a FODC over $\A $, $\sigma $ a braiding
and $g$ a $\sigma $-metric of $\Gam $.
Similarly to the construction in \cite{a-Heck99} one can define a contraction
mapping $\ctr{\cdot }{\cdot } :\Gamt[k] \otA \Gamt[l] \to
\Gamt[\mathrm{max}\{l-k,0\}]$. Differing from the notation therein we set
$\ctr{\rho }{\rho '}=0$ for $\rho \in \Gamt[k]$, $\rho '\in \Gamt[l]$, $k>l$,
and retain the formulas
\begin{align}
\ctr{\rho }{\rho '}&=g_{12}(\rho \otA B_{1,l-1}(\rho ')),&
\ctr{\rho \otA \rho ''}{\rho '}&=\ctr{\rho }{\ctr{\rho ''}{\rho '}\,},
\end{align}
where $\rho \in \Gam $, $\rho '\in \Gamt[l]$, $\rho ''\in \Gamt $.
In \cite{a-BCDRV96} Bautista et.\,al.~proposed the notion of a quantum
Clifford algebra. Adapting their ideas to the present situation and having
the classical situation in mind we introduce the following definition.

\begin{defin}\label{d-qclifalg}
Let $\mc{I}$ denote the two-sided ideal of the tensor algebra $\Gamt $
generated by the elements
\begin{align}\label{eq-clifideal}
\{\rho _2-g(\rho _2)\,|\,\rho _2\in \ker (\id -\sigma )\}.
\end{align}
The algebra $\Gamt /\mc{I}$ is called
the \textit{quantum Clifford algebra for the left-cova\-riant $\A $-bimodule
$\Gam $ and the $\sigma $-metric $g$}. We denote it by
$\Cl (\Gam ,\sigma ,g)$.
\end{defin}

Let $\multcl $ denote the canonical mapping
$\multcl :\Gamt[2]/(\mc{I}\cap (\A \oplus \Gamt[2]))$.
Similarly to the classical case, $\Cl (\Gam ,\sigma ,g)$ can be equipped with
a $\mathbb{Z}_2$-grading such that $\mathrm{deg}\,\rho _k=(-1)^k$ for
$\rho _k\in \Gamt[k]$.
Let $\Cl ^\eta (\Gam ,\sigma ,g)$, $\eta \in \{+,-\}$, denote the subspace
of homogeneous elements of degree $\eta $.

Since $g$ is left-covariant, there exists a natural left coaction of $\A $
on $\Cl (\Gam ,\sigma ,g)$.

In the classical situation, the Clifford algebra has a representation over
the exterior algebra. Now we should prove that this is valid also in the
quantum case.

\begin{satz}
The formulas
\begin{align}\label{eq-clifrep}
a\lact \rho &:=a\rho,&
\omega \lact \rho &:=\omega \wedge \rho +\ctr{\omega }{\rho },
\end{align}
where $a\in \A $, $\rho \in \Gamw $ and $\omega \in \Gam $,
define a representation of $\Cl (\Gam ,\sigma ,g)$ over $\Gamw $.
\end{satz}

\begin{bew}
Since $\ctr{\cdot }{\cdot }$ is a homomorphism of the left $\A $-modules
$\Gam \otA \Gamw $ and $\Gamw $, it is easy to see that
(\ref{eq-clifrep}) gives a well-defined action of the
$\A $-bimodule $\Gam $ on $\Gamw $. Hence this action can be extended
to an action $\lact $ on $\Gamw $ of the algebra $\Gamt $. Now we only
have to prove that the elements in (\ref{eq-clifideal}) act trivially
on $\Gamw $.

Suppose that $a_{ij}\w _i\otA \w _j\in \ker (\id -\sigma )$ and
$\rho _k\in \Gamw[k]$, $k\geq 0$. Then
\begin{align*}
&(a_{ij}\w _i\otA \w _j-a_{ij}g_{ij})\lact \rho _k
=a_{ij}\w _i\lact (\w _j\lact \rho _k)-a_{ij}g_{ij}\rho _k\\
&\quad =a_{ij}\w _i\lact \left(\w _j\wedge \rho _k+\ctr{\w _j}{\rho _k}\right)
 -a_{ij}g_{ij}\rho _k\\
&\quad =a_{ij}\w _i\wedge \w _j\wedge \rho _k
 +\ctr{a_{ij}\w _i}{\w _j\wedge \rho _k}\\
&\quad \phantom{=}+a_{ij}\w _i\wedge \ctr{\w _j}{\rho _k}
 +\ctr{a_{ij}\w _i}{\ctr{\w _j}{\rho _k}\,}
 -a_{ij}g_{ij}\rho _k.\tag{$*$}
\end{align*}
Since $\ctr{a_{ij}\w _i}{\ctr{\w _j}{\rho _k}\,}=
\ctr{a_{ij}\w _i\wedge \w _j}{\rho _k}$, the first and fourth summands
of the above expression vanish.
By formula (39) in \cite{a-Heck99} the second summand can be written as
\begin{align*}
a_{ij}\ctr{\w _i}{\w _j}\rho _k-b_{rs}\w _r\wedge \ctr{\w _s}{\rho _k},
\tag{$**$}
\end{align*}
where $b_{rs}\w _r\otA \w _s=\sigma ^{-1}(a_{ij}\w _i\otA \w _j)$.
The first summand of ($**$) is the same as the last one in ($*$).
Therefore ($*$) becomes 
\begin{align*}
-b_{rs}\w _r\wedge \ctr{\w _s}{\rho _k}
+a_{ij}\w _i\wedge \ctr{\w _j}{\rho _k}
=(a_{ij}-b_{ij})\w _i\wedge \ctr{\w _j}{\rho _k}.
\end{align*}
But $(\id -\sigma ^{-1})(a_{ij}\w _i\otA \w _j)=
-\sigma ^{-1}(\id -\sigma )(a_{ij}\w _i\otA \w _j)=0$, and so
$a_{ij}=b_{ij}$ for any $i,j=1,\ldots ,n$. This means that ($*$) vanishes
for any $\rho _k\in \Gamw[k]$.
\end{bew}

\begin{bem}
Since $\rho \lact 1=\rho $ for any $\rho \in \A \oplus \Gam $,
$\rho =0$ in $\Cl (\Gam ,\sigma ,g)$ implies that $\rho =0$ in
$\A \oplus \Gam $.
Hence $\A $ and $\Gam $ can be naturaly embedded in $\Cl (\Gam ,\sigma ,g)$.
\end{bem}

\begin{satz}\label{s-betaCl}
The involution $\beta $ of the $\A $-bimodule $\Gam $ extends uniquely to
an involution of the quantum Clifford algebra $\Cl (\Gam ,\sigma ,g)$.
\end{satz}

\begin{bew}
We have to show that if $(\id -\sigma )\rho _2=0$ for
$\rho _2=a_{ij}\w _i\otA \w _j\in \Gamt[2]$
then $\beta (\w _j)\beta (\w _i)\beta (a_{ij})=\beta (g(\rho _2))$
in $\Cl (\Gam ,\sigma ,g)$. Since $\ker (\id -\sigma )$ is a left-covariant
left $\A $-module, we can assume that $a_{ij}\in \comp $.
By definition we have
\begin{align*}
\beta (\w _j)\otA \beta (\w _i)&=
\overline{B^k_j}\overline{B^l_i}\w _k\otA \w _l=
F^k_m\overline{E^m_j}F^l_n\overline{E^n_i}\w _k\otA \w _l\\
&=\overline{E^m_j}\overline{E^n_i}S^2(\w _m\otA \w _n)
=S^2(\w _j^*\otA \w _i^*).
\end{align*}
Since $S^2\sigma =\sigma S^2$ and
$\sigma (\rho ^*)=(\sigma ^\eta (\rho ))^*$ for $\rho \in \Gamt[2]$
with $\eta =+1$ or $\eta =-1$,
it follows that
$(\id -\sigma )(\beta (\w _j)\otA \beta (\w _i)\overline{a_{ij}})=
\beta (\id -\sigma ^\eta )(\rho _2)=0$.
This means that $\beta (\w _j)\beta (\w _i)\overline{a_{ij}}=
g(\beta (\w _j)\otA \beta (\w _i)\overline{a_{ij}})$ in
$\Cl (\Gam ,\sigma ,g)$. But the latter is equal to $\beta (g(\rho _2))$
because of $gS^2=S^2g$ and the requirement (\ref{eq-gstar}).
\end{bew}

Similarly to the classical situation a left-covariant left ideal
$\mc{S}$ of $\Cl (\Gam ,\sigma ,g)$ is called a \textit{spinor module}
if $\mc{S}$ is an $\A $-bimodule, too. By Theorem 4.1.1.\ in
\cite{b-Sweedler} it follows that $\mc{S}$ is a free left
$\A $-module and there exists a left-invariant $\A $-module basis of $\mc{S}$.
Let $\multcls $ denote the left action
$\multcls :\Gam \otA \mc{S}\to \mc{S}$ of $\Gam $ on the spinor module
$\mc{S}$.

One of the most important aims of this article is to introduce a
(hermitean non-degenerate) metric on spinor modules $\mc{S}$.
In order to do this let
$\cconj{\mc{S}}$ denote the complex conjugate of the vector space $\mc{S}$,
that is $\cconj{\mc{S}}=\mc{S}$ as sets and
$\lambda \cconj{\psi }=\cconj{(\bar{\lambda }\psi )}$ for
$\lambda \in \comp $, $\psi \in \mc{S}$.
We define a left coaction $\lkow $ and an $\A $-bimodule structure on the
vector space $\mc{S}\ot \cconj{\mc{S}}$ by
\begin{equation}
\begin{gathered}
\lkow (\psi \ot \cconj{\psi '{}}):=
\psi _{(-1)}(\psi '_{(-1)})^*\ot (\psi _{(0)}\ot \cconj{\psi '_{(0)}{}}),\\
a(\psi \ot \cconj{\psi '{}})=a\psi \ot \cconj{\psi '{}},\qquad
(\psi \ot \cconj{\psi '{}})a=\psi \ot \cconj{a^*\psi '{}}
\end{gathered}
\end{equation}
for $\psi ,\psi '\in \mc{S}$, $a\in \A $,
where $\lkow (\psi )=\psi _{(-1)}\ot \psi _{(0)}$,
$\lkow (\psi ')=\psi '_{(-1)}\ot \psi '_{(0)}$ and $*$ is the involution
of $\A $.

\begin{defin}\label{d-hmet}
A mapping $\hmet{\cdot }{\cdot }:\mc{S}\ot \cconj{\mc{S}}\to \comp $ is called
a \textit{metric on} $\mc{S}$ if\\
\indent (i) there exists a left-covariant mapping
$\hmet{\cdot }{\cdot }_0:\mc{S}\ot \cconj{\mc{S}}\to \A $ such that both
$\hmet{a\psi }{\cconj{(b\psi ')}}_0=a\hmet{\psi }{\cconj{\psi '{}}}_0b^*$
for any $a,b\in \A $ and $\psi ,\psi '\in \mc{S}$ and the equation
$\hmet{\cdot }{\cdot }=h\circ \hmet{\cdot }{\cdot }_0$ hold,\\
\indent (ii) $\hmet{\cdot }{\cdot}$ is non-degenerate and hermitean
(but it needs not be positive definite).\\
To simplify the notation we will write $\hmet{a\psi }{b\psi '}$
for $\hmet{a\psi }{\cconj{(b\psi ')}}$. But remember that the metric is
antilinear in the second component.
\end{defin}

\begin{bems}
1. Let $\hmet{\cdot }{\cdot }_0:\mc{S}\ot \cconj{\mc{S}}\to \A $ be
such that $\hmet{a\psi }{b\psi '}_0=a\hmet{\psi }{\psi '}_0b^*$
for any $a,b\in \A $ and $\psi ,\psi '\in \mc{S}$.
Then it is left-covariant if and only if
$\hmet{\psi }{\psi '}_0\in \comp $
for any $\psi ,\psi '\in \linv{\mc{S}}$.
Hence, if $\hmet{\cdot }{\cdot }$ is a metric on $\mc{S}$ then
$\hmet{\cdot }{\cdot }_0$ can be reconstructed by the setting
$\hmet{a\psi }{b\psi '}_0:=a\hmet{\psi }{\psi '}b^*$ for $a,b\in \A $,
$\psi ,\psi '\in \linv{\mc{S}}$.
\\
2. A mapping $\hmet{\cdot }{\cdot }$ in (i) of Definition \ref{d-hmet} is
non-degenerate if and only if its restriction
$\hmet{\cdot }{\cdot }:\linv{\mc{S}}\ot \cconj{\linv{\mc{S}}}\to \comp $
is non-degenerate and\\
3. it is hermitean if and only if its restriction
$\hmet{\cdot }{\cdot }:\linv{\mc{S}}\ot \cconj{\linv{\mc{S}}}\to \comp $
is hermitean.\\
4. Any metric $\hmet{\cdot }{\cdot }$ on $\mc{S}$ satisfies the equations
\begin{align}
\hmet{a\psi }{b\psi '}_0^*&=\hmet{b\psi '}{a\psi }_0\quad \text{and}\\
\hmet{a\psi }{\psi '}&=\hmet{\psi }{\beta (a)\psi '}
\end{align}
for all $a\in \A $ and $\psi ,\psi '\in \mc{S}$.
\end{bems}

Of course, metrics on spinor modules always exist. In this paper interesting
metrics will also have to satisfy the equation
\begin{align}\label{eq-hmetbeta}
\hmet{u\varphi }{\psi }&=\hmet{\varphi }{\beta (u)\psi }
\end{align}
for $u\in \Cl (\Gam ,\sigma ,g)$ and $\varphi ,\psi \in \mc{S}$
(see Proposition \ref{s-symmdirac}).
A similar construction as in Section 8.2.1 in \cite{b-Crumey90} gives
the existence of such a metric on ''minimal'' spinor modules
(minimal left ideals).

Observe that $\beta (\mc{S})$ is a minimal right ideal.
If $\Cl (\Gam ,\sigma ,g)$ is
a simple left-covariant algebra, then $\beta (\linv{\mc{S}})\linv{\mc{S}}=
\beta (\linv{\mc{S}})\cap \linv{\mc{S}}$
is a one-dimensional subspace of $\linv{\mc{S}}$. Since $\beta ^2=\id $,
it contains an element $x$ such that $\beta (x)=x$.
Consider elements $\psi ,\psi '\in \mc{S}$. Since $\mc{S}$ is an
$\A $-bimodule, there exist $a,b\in \A $ such that $\beta (\psi ')\psi =ax=xb$.
Hence $xc=c_{(1)}\hat{f}(c_{(2)})x$ for all $c\in \A $ and therefore
$h(a)=h(b)$. We define $\hmet{\psi }{\psi '}=h(a)=h(b)$.
Since
\begin{align*}
\beta (\psi )\psi '&=\beta (\beta (\psi ')\psi )=\beta (ax)=x\beta (a),
\end{align*}
we obtain $\hmet{\psi '}{\psi }=h(\beta (a))=h(\rho (a)^*)=\overline{h(a)}$.
This means that the mapping $\hmet{\cdot }{\cdot }$ is hermitean.
For an arbitrary element $u\in \Cl (\Gam ,\sigma ,g)$ we get
\begin{align*}
\beta (\psi ')u\psi &=\beta (\beta (u)\psi ')\psi 
\end{align*}
and therefore $\hmet{u\psi }{\psi '}=\hmet{\psi }{\beta (u)\psi '}$.
Finally we have to prove the non-degeneracy of $\hmet{\cdot }{\cdot }$.
Since $\beta ({\linv{\mc{S}}})\linv{\mc{S}}$ is a nontrivial vector space,
there are $\psi ',\psi ''\in \linv{\mc{S}}$ such that
$\beta (\psi '')\psi '\not= 0$. Because $\mc{S}$ is a minimal left ideal
of $\Cl (\Gam ,\sigma ,g)$, for any $\psi \in \linv{\mc{S}}$, $\psi \not=0$,
there exists $u\in \Cl (\Gam ,\sigma ,g)$ such that $\psi ''=u\psi $.
Then $0\not=\hmet{\psi '}{\psi ''}=\hmet{\psi '}{u\psi }
=\hmet{\beta (u)\psi '}{\psi }$ and hence $\hmet{\cdot }{\cdot }$ is
non-degenerate on $\linv{\mc{S}}\ot \linv{\cconj{\mc{S}}}$.
This means that there are bases $\{\psi _i\}$ and $\{\psi '_i\}$ of
$\linv{\mc{S}}$ such that $\hmet{\psi _i}{\psi '_j}=\delta _{ij}$.

Suppose now that $\psi \in \mc{S}\setminus \{0\}$. Then $\psi =\psi _ia_i$,
$a_i\in \A $, and we can assume that $a:=a_1\not=0$. There exists $b\in \A $
such that $h(ba_{(1)}\hat{f}(a_{(2)}))\not=0$. We conclude that
\begin{align*}
\beta (\psi '_1\beta (b))\psi &=b\beta (\psi '_1)\psi _ia_i=bxa
=ba_{(1)}\hat{f}(a_{(2)})x
\end{align*}
from which $\hmet{\psi }{\psi '_1\beta (b)}\not=0$ follows.
This proves the following proposition.

\begin{satz}\label{s-spinmet}
Suppose that $\Cl (\Gam ,\sigma ,g)$ is a simple left-covariant algebra
and $\mc{S}$ is a minimal left ideal of $\Cl (\Gam ,\sigma ,g)$ and an
$\A $-bimodule. Then there exists a metric on $\mc{S}$ such that
$\hmet{u\psi }{\psi '}=\hmet{\psi }{\beta (u)\psi '}$ holds for
$u\in \Cl (\Gam ,\sigma ,g)$, $\psi ,\psi '\in \mc{S}$.
\end{satz}

\section{Connections}
\label{sec-connection}

Let $\mc{C}$ be an $\A $-bimodule and $(\Gam ,\dif )$ a first order
differential calculus over $\A $.
Following \cite{a-BDMS1} we call a map
$\nabla :\mc{C} \to \Gam \otA \mc{C}$
a \textit{left connection on} $\mc{C}$ if it satisfies the rule
\begin{align}\label{eq-connection}
\nabla (a\rho )&=\dif a \otA \rho +a\nabla (\rho )
\end{align}
for any $a\in \A $ and any $\rho \in \mc{C}$.
If $\mc{C}$ is a left-covariant $\A $-bimodule then a left connection
$\nabla $ on $\mc{C}$ is called \textit{left-covariant},
if $\lkow (\nabla (\rho ))=(\id \ot \nabla )\lkow (\rho )$.

Unfortunately, in general it is
not possible to extend a connection on an $\A $-bimodule $\mc{C}$ to the
tensor product $\mc{C}\otA \mc{C}$ ($\mc{C}^\otA $, respectively).
In \cite{a-BDMS1} a constraint was
given in which case such an extension can be made. The definition therein
yields that a left connection is extensible if and only if there exists a
bimodule homomorphism $\ts :\mc{C}\otA \Gam \to \Gam \otA \mc{C}$ such that
\begin{align}
\ts (\rho \otA \dif a)=\nabla (\rho a)-\nabla (\rho )a
\end{align}
holds for all $a\in \A $ and $\rho \in \mc{C}$.
In this case, following Mourad's definition \cite{a-Mourad1,a-DMMM95},
$\nabla $ is called a \textit{linear left connection on} the $\A $-bimodule
$\mc{C}$.

Let $\nabla $ be a linear left connection on $\mc{C}$. One defines the
left connection
$\nabla :\mc{C}^{\otA k}\to \Gam \otA \mc{C}^{\otA k}$, $k\geq 2$,
recursively by
\begin{align}\label{eq-extconn}
\nabla (\rho \otA \rho _{k-1}):=\nabla (\rho )\otA \rho _{k-1}
+(\ts \otA \id ^{k-1})(\rho \otA \nabla (\rho _{k-1}))
\end{align}
for $\rho \in \mc{C}$ and $\rho _{k-1}\in \mc{C}^{\otA k-1}$.
Moreover, if $\mc{C}'$ is another $\A $-bimodule and $\nabla '$ is a left
connection on $\mc{C}'$ then the formula
\begin{align}\label{eq-tensorconn}
\widetilde{\nabla }(\rho \otA \rho '):=\nabla (\rho )\otA \rho '
+(\ts \otA \id )(\rho \otA \nabla '(\rho ')),
\end{align}
$\rho \in \mc{C}$, $\rho '\in \mc{C}'$, defines a left connection on
$\mc{C}\otA \mc{C}'$. If additionally $\nabla '$ is an extensible left
connection on $\mc{C}'$,
$\nabla '(\rho 'a)=\nabla '(\rho ')a+\ts '(\rho '\otA \dif a)$, then
$\widetilde{\nabla }$ satisfies the equation
\begin{align}
\widetilde{\nabla }(\rho \otA \rho 'a)=\widetilde{\nabla }(\rho \otA \rho ')a
+(\ts \ot \id)(\id \ot \ts ')(\rho \otA \rho '\otA \dif a)
\end{align}
for $\rho \in \mc{C},\rho '\in \mc{C}'$ and $a\in \A $.

Let $\mc{C}$ be a left-covariant $\A $-bimodule and let
$\{\cb _\alpha \,|\,\alpha =1,\ldots ,M:=\dim \mc{C}\}$ be a basis of
$\linv{\mc{C}}:=\{\cb \in \mc{C}\,|\,\lkow (\cb )=1\ot \cb \}$.
There are linear functionals $F^\alpha _\beta $ on $\A $ such that
$\cb _\alpha a=a_{(1)}F^\alpha _\beta (a_{(2)})\cb _\beta $ for any
$a\in \A $, $\alpha =1,\ldots ,M$.

\begin{lemma}\label{l-extlconn}
Let $\nabla $ be a left-covariant left connection on $\mc{C}$,
$\nabla (\cb _\alpha )=\crs ^{i\beta }_\alpha \w _i\otA \cb _\beta $,
$\crs ^{i\beta }_\alpha \in \comp $.
Then $\nabla $ is a linear left connection on $\mc{C}$ if and only if
\begin{align}
\crs ^{i\beta }_\alpha \koun +S(F^\alpha _\gamma )\x _iF^\gamma _\beta
-S(F^\alpha _\delta )\crs ^{j\gamma }_\delta f^j_i F^\gamma _\beta \in \mc{X}
\end{align}
for any $\alpha ,\beta =1,\ldots ,M$, $i=1,\ldots ,n$.
\end{lemma}

\begin{bew}
Recall that $\{\w _i\otA \cb _\beta \}$ and $\{\cb _\alpha \}$ are free
bases of the left $\A $-modules $\Gam \otA \mc{C}$ and $\mc{C}$, respectively.
Hence by (\ref{eq-connection}) any left connection on $\mc{C}$ is uniquely
determined by elements $\crs ^{i\beta }_\alpha \in \A $.
Moreover, left-covariance of $\nabla $ is satisfied if and only if
$\crs ^{i\beta }_\alpha \in \comp $ for any
$\alpha ,\beta =1,\ldots ,\dim \mc{C}$ and $i=1,\ldots ,n$.

By the above considerations and because of
$\dif a=a_{(1)}S(a_{(2)})\dif a_{(3)}$, $\nabla $ is a linear left connection
if and only if the mapping $\ts $, given by
\begin{equation}
\begin{aligned}
\ts (\rho _iS(a_i{}_{(1)})\otA \dif a_i{}_{(2)}):=&
\nabla (\rho _iS(a_i{}_{(1)})a_i{}_{(2)})
-\nabla (\rho _iS(a_i{}_{(1)}))a_i{}_{(2)}\\
=&\koun (a_i)\nabla (\rho _i)-\nabla (\rho _iS(a_i{}_{(1)}))a_i{}_{(2)},
\end{aligned}
\end{equation}
$\rho _i\in \mc{C}$, $a_i\in \A $, is well defined.
On the left hand side the expression
$\ts (\rho _i\otA S(a_i{}_{(1)})\dif a_i{}_{(2)})=
\ts (\rho _i\otA \omega (a_i))$
appears which is zero for linearly independent $\rho _i$ if and only if
$a_i\in \mc{R}_\Gam +\comp \cdot 1$.
This gives that $\ts $ is well-defined if and only if
$\koun (a)\nabla (\rho )-\nabla (\rho S(a_{(1)}))a_{(2)}=0$
for any $a\in \mc{R}_\Gam +\comp \cdot 1$ and $\rho \in \linv{\mc{C}}$.
Now let us compute this expression for $\rho =\cb _\alpha $.
\begin{align*}
&\koun (a)\nabla (\cb _\alpha )-\nabla (\cb _\alpha S(a_{(1)}))a_{(2)}=\\
&\quad
=\koun (a)\crs ^{i\beta }_\alpha \w _i\otA \cb _\beta
-\nabla (S(a_{(2)})F^\alpha _\gamma (S(a_{(1)}))\cb _\gamma )a_{(3)}\\
&\quad \tag{$*$}
=\koun (a)\crs ^{i\beta }_\alpha \w _i\otA \cb _\beta
-\dif S(a_{(2)})S(F^\alpha _\gamma )(a_{(1)})\otA \cb _\gamma a_{(3)}\\
&\quad \phantom{=}
-S(a_{(2)})S(F^\alpha _\delta )(a_{(1)})\crs ^{j\gamma }_\delta
 \w _j\otA \cb _\gamma a_{(3)}.
\end{align*}
The second summand of ($*$) can be reformulated as
\begin{align*}
&-\dif S(a_{(2)})S(F^\alpha _\gamma )(a_{(1)})\otA a_{(3)}
 F^\gamma _\beta (a_{(4)})\cb _\beta =\\
&\qquad
=(-\dif (S(a_{(2)})a_{(3)})+S(a_{(2)})\dif a_{(3)})
 S(F^\alpha _\gamma )(a_{(1)})F^\gamma _\beta (a_{(4)})\otA \cb _\beta \\
&\qquad
=S(a_{(2)})a_{(3)}\x _i(a_{(4)})\w _i
 S(F^\alpha _\gamma )(a_{(1)})F^\gamma _\beta (a_{(5)})\otA \cb _\beta \\
&\qquad
=(S(F^\alpha _\gamma )\x _iF^\gamma _\beta )(a)\w _i \otA \cb _\beta .
\end{align*}
For the third summand of ($*$) one computes
\begin{align*}
&-\crs ^{j\gamma }_\delta S(a_{(2)})S(F^\alpha _\delta )(a_{(1)})a_{(3)}
 (f^j_i F^\gamma _\beta )(a_{(4)}) \w _i\otA \cb _\beta =\\
&\qquad
=-\crs ^{j\gamma }_\delta (S(F^\alpha _\delta )f^j_i F^\gamma _\beta )(a)
 \w _i\otA \cb _\beta .
\end{align*}
Together we obtain
\begin{multline}
\koun (a)\nabla (\cb _\alpha )-\nabla (\cb _\alpha S(a_{(1)}))a_{(2)}=\\
=(\crs ^{i\beta }_\alpha \koun +S(F^\alpha _\gamma )\x _iF^\gamma _\beta
-S(F^\alpha _\delta )\crs ^{j\gamma }_\delta f^j_i F^\gamma _\beta )(a)
\w _i\otA \cb _\beta
\end{multline}
for any $a\in \A $ and $\alpha =1,\ldots ,M$.
Finally, recall that $f\in \A '$, $f(a)=0$ for any $a\in \mc{R}_\Gam
+\comp \cdot 1$ is by definition equivalent to $f\in \mc{X}$.
\end{bew}

\begin{defin}
Let $\mc{C}$ be an $\A $-bimodule algebra with multiplication $\circ $.
Then a linear mapping
$\nabla :\mc{C}\to \Gam \ot \mc{C}$ is called a
\textit{linear left connection on}
$\mc{C}$, if\\
\indent (i) $\nabla $ is a left connection on the $\A $-bimodule
$\mc{C}$,\\
\indent (ii) there exists a bimodule homomorphism
$\ts :\mc{C}\ot \Gam \to \Gam \ot \mc{C}$ such that
$\ts (\rho \otA \dif a)=\nabla (\rho a)-\nabla (\rho )a$
for any $a\in \A $, $\rho \in \mc{C}$, and\\
\indent (iii) $\nabla (\rho \circ \rho ')=
(\id \otA \circ )\bigl( \nabla (\rho )\otA \rho '
+(\ts \otA \id )(\rho \otA \nabla (\rho '))\bigr) $
holds for any $\rho ,\rho '\in \mc{C}$.
\end{defin}

Let $\mc{C}$ be an $\A $-bimodule algebra and $\mc{S}$ a left $\mc{C}$-module
with left action $\lact $.
Then a mapping $\szushg :\mc{S}\to \Gam \otA \mc{S}$ is called a
\textit{linear left connection on} $\mc{S}$ if there is a left connection
$\nabla $ on $\mc{C}$ such that
\begin{align}\label{eq-szushg}
\szushg (a\varphi )&=(\id \otA \lact )(\nabla (a)\otA \varphi
+(\ts \otA \id )(a\otA \szushg (\varphi )))
\end{align}
for any $a\in \mc{C}$ and $\varphi \in \mc{S}$.

Let us call a FODC $(\Gam ,\dif )$ over $\A $ \textit{inner},
if there exists an $\omega \in \Gam $ such that
$\dif a=\omega a-a\omega $ for any $a\in \A $.

\begin{lemma}\label{l-innzushg}
Let $\mc{C}$ be an $\A $-bimodule and
let $(\Gam ,\dif )$ be an inner FODC over $\A $
such that $\dif a=\omega a-a\omega $ for $a\in \A $.
Let $\ts :\mc{C} \otA \Gam \to \Gam \otA \mc{C} $ and
$V:\mc{C} \to \Gam \otA \mc{C} $ be homomorphisms of $\A $-bimodules.
Then the assignment
\begin{align}\label{eq-innzushg}
\nabla (\rho ):=\omega \otA \rho -\ts (\rho \otA \omega )
+V(\rho )
\end{align}
defines a linear left connection $\nabla $ on $\mc{C} $ such that
$\nabla (\rho a)-\nabla (\rho )a =\ts (\rho \otA \dif a)$.
Moreover, any linear left connection on $\mc{C}$ is given in this manner.
\end{lemma}

\begin{bem}
For $\mc{C}=\Gam $ the assertion of the lemma was already proved in
\cite[Prop.~C.3]{a-BDMS1}. Since the proof of this lemma is similar,
we omit it.
\end{bem}

For inner first order differential calculi $\Gam $ over $\A $
the following proposition holds.

\begin{satz}\label{s-clconn}
Suppose that $(\Gam ,\dif )$ is an inner FODC and $\nabla $ is a linear
left connection on $\Gam $ such that
$\nabla (\rho )=\omega \otA \rho -\ts (\rho \otA \omega )$,
$\rho \in \Gam $, for an endomorphism $\ts $ of the $\A $-bimodule
$\Gam \otA \Gam $.\\
1.\ The linear left connection $\nabla :\Gamt \to \Gam \otA \Gamt $ is
compatible with Wo\-ro\-no\-wicz' antisymmetrizer if
$\sigma _{23}\ts _{12}\ts _{23}=
\ts _{12}\ts _{23}\sigma _{12}$.\\
2.\ The linear left connection $\nabla :\Gamt \to \Gam \otA \Gamt $ is
compatible with the multiplication of the Clifford algebra if
$\sigma _{23}\ts _{12}\ts _{23}=\ts _{12}\ts _{23}\sigma _{12}$ and
$g_{23}\ts _{12}\ts _{23}=g_{12}$.
\end{satz}

\begin{bew}
First one checks with (\ref{eq-extconn}) that
\begin{align}
\nabla (\rho _k)=\omega \otA \rho _k-\ts _{12}\ts _{23}
\cdots \ts _{k,k+1}(\rho _k\otA \omega )
\end{align}
for any $\rho _k\in \Gamt[k]$. Let
$\ts :\Gamt\otA \Gam \to \Gam \otA \Gamt $ denote the linear mapping
defined by
\begin{align}
\ts (\rho _k\otA \rho ):=\ts _{12}\ts _{23}\cdots \ts _{k,k+1}
(\rho _k\otA \rho ),\quad \rho _k\in \Gamt[k],\rho \in \Gam .
\end{align}

Suppose that $\sigma _{23}\ts _{12}\ts _{23}=\ts _{12}\ts _{23}\sigma _{12}$.
Then for any $i\in \mathbb{N}$, $1<i\leq k$ we have
$\sigma _{i,i+1}\ts _{12}\ts _{23}\cdots \ts _{k,k+1}=
\ts _{12}\ts _{23}\cdots \ts _{k,k+1}\sigma _{i-1,i}$
and hence
\begin{align}
(\id \otA A_k)\ts _{12}\ts _{23}\cdots \ts _{k,k+1}=
\ts _{12}\ts _{23}\cdots \ts _{k,k+1}(A_k\otA \id )
\quad \text{for any $k\geq 2$}.
\end{align}
Similarly, if $g_{23}\ts _{12}\ts _{23}=g_{12}$ and $k\geq 2$ then
for any $i\in \mathbb{N}$, $1<i\leq k$ we have
$g_{i,i+1}\ts _{12}\cdots \ts _{k,k+1}=
\ts _{12}\cdots \ts _{k-2,k-1}g_{i-1,i}$.

Suppose now that $A_k(\rho _k)=0$ for a $\rho _k\in \Gamt[k]$, $k\geq 2$.
Then
\begin{align*}
(\id \otA A_k)\nabla (\rho _k)=&
\omega \otA A_k(\rho _k)-(\id \otA A_k)\ts (\rho _k\otA \omega )\\
=&-\ts (A_k(\rho _k)\otA \omega )=0.
\end{align*}
Hence $\nabla $ is a well-defined connection on $\Gamw $.

Let now $\rho \in \Gamt $ be an element of the ideal $\mc{I}$ in Definition
\ref{d-qclifalg}. Without loss of generality we may assume that
$\rho =\rho '\otA (\rho _2-g(\rho _2))\otA \rho ''$, where
$\rho '\in \Gamt[k]$, $\rho ''\in \Gamt[l]$, $k,l\in \mathbb{N}_0$ and
$\rho _2\in \Gamt[2]$, $\sigma (\rho _2)=\rho _2$.
Then $\nabla (\rho )=\omega \otA \rho -\ts (\rho \otA \omega )$ and the first
summand is an element of $\Gam \otA \mc{I}$. Now it suffices to show that
\begin{gather}\label{eq-rhoinI1}
(\id -\sigma _{k+2,k+3})\ts ((\rho '\otA \rho _2\otA \rho '')\otA \omega )=0
\qquad \text{and}\\
\label{eq-rhoinI2}
g_{k+2,k+3}\ts ((\rho '\otA \rho _2\otA \rho '')\otA \omega )=
\ts ((\rho 'g(\rho _2)\otA \rho '')\otA \omega ).
\end{gather}
But (\ref{eq-rhoinI1}) follows from
$\sigma _{k+2,k+3}\ts =\ts \sigma _{k+1,k+2}$ and $\sigma \rho _2=\rho _2$
and (\ref{eq-rhoinI2}) is proved by $g_{k+2,k+3}\ts =\ts g_{k+1,k+2}$.
\end{bew}

\begin{satz}\label{s-modcompconn}
Let $(\Gam ,\dif )$ be an inner FODC over $\A $, $\dif a=\omega a-a\omega $
for $a\in \A $. Let $\mc{C}$ be an
$\A $-bimodule and left $\Gamt $-module with left action $\mathrm{m}$.
Suppose that $\nabla $ and $\nabla '$ are linear left connections on $\Gamt $
and on the $\A $-bimodule $\mc{C}$, respectively, such that
\begin{align}
\nabla (\rho )&=\omega \otA \rho -\ts (\rho \otA \omega ),&
\nabla '(\psi )&=\omega \otA \psi -\tau (\psi \otA \omega )+V(\psi )
\end{align}
holds for $\rho \in \Gamt ,\psi \in \mc{C}$. Then $\nabla '$ is
compatible with the left action $\mathrm{m}$ of $\Gamt $ on $\mc{C}$
if and only if
\begin{align}\label{eq-modcompconn}
(\tau \mathrm{m}_{12}-\mathrm{m}_{23}\ts _{12}\tau _{23})
(\rho \otA \psi \otA \omega )
+(\mathrm{m}_{23}\ts _{12}V_2-V\mathrm{m})(\rho \otA \psi )&=0
\end{align}
for all $\rho \in \Gam ,\psi \in \mc{C}$.
\end{satz}

\begin{bew}
It is enough to compare both sides of (\ref{eq-szushg}) for $a\in \Gam $
and $\varphi \in \mc{C}$.
The left hand side becomes $\omega \otA a\varphi -\tau (a\varphi \otA \omega )
+V(a\varphi )$. The right hand side takes the form
\begin{align*}
&\mathrm{m}_{23}(\nabla (a)\otA \varphi +\ts _{12}(a\otA \nabla '\varphi ))=\\
&\quad =\mathrm{m}_{23}(\omega \otA a\otA \varphi
-\ts _{12}(a\otA \omega \otA \varphi )\\
&\quad \phantom{=\mathrm{m}_{23}()}
+\ts _{12}(a\otA \omega \otA \varphi )
-\ts _{12}\tau _{23}(a\otA \varphi \otA \omega )
+\ts _{12}V_2(a\otA \varphi ))\\
&\quad =\omega \otA a\varphi
-\mathrm{m}_{23}\ts _{12}\tau _{23}(a\otA \varphi \otA \omega )
+\mathrm{m}_{23}\ts _{12}V_2(a\otA \varphi ).
\end{align*}
Both sides are equal if and only if (\ref{eq-modcompconn}) holds.
\end{bew}

Let $\Gamw $ be a differential calculus over $\A $ with first order part
$\Gam $. Let $\nabla $ be a left connection on $\Gam $. The mapping
$\tors :\Gam \to \Gamw[2]$, defined by the formula
$\tors :=\multw \nabla -\dif $, is called the \textit{torsion} of $\nabla $.
It satisfies $\tors (a\rho )=a\tors (\rho )$
for all $a\in \A ,\rho \in \Gam $. If $\nabla $ is a linear left connection
on $\Gam $, $\nabla (\rho a)=\nabla (\rho )a+\ts (\rho \otA \dif a)$,
then $\tors $ is a homomorphism of $\A $-bimodules if and only if
\begin{align}\label{eq-torshom}
\multw (\id +\ts )(\rho \otA \rho ')&=0
\end{align}
for all $\rho ,\rho '\in \Gam $.

Suppose that $\Gamw $ is an inner differential calculus over $\A $, that is
there exists $\omega \in \Gam $ such that $\dif \rho =\omega \wedge \rho
-(-1)^k\rho \wedge \omega $ for all $\rho \in \Gamw[k]$.
Let $\nabla $ be a linear left connection on $\Gam $ which satisfies
$\nabla \rho =\omega \otA \rho -\ts (\rho \otA \omega )$
for all $\rho \in \Gam $. If the torsion $\tors $ of $\nabla $ fulfills
(\ref{eq-torshom}) then $\tors (\rho )=0$ for all $\rho \in \Gam $.

Let $\mc{S}$ be a spinor module, $\hmet{\cdot }{\cdot }$ a metric on $\mc{S}$
and let $\nabla $ and $\szushg $ be linear left connections on
$\Cl (\Gam ,\sigma, g)$ and on $\mc{S}$, respectively.
Generalizing the notion of Definition \ref{d-hmet}, we use the symbol
$\hmet{a\psi }{b\psi '}_0:=a\otA \hmet{\psi }{\psi '}_0b^*$
for all $a,b\in \Cl (\Gam ,\sigma ,g)$, $\psi ,\psi '\in \mc{S}$.
Then the mapping $\szushg ^*: \mc{S}\to \Gam \otA \mc{S}$, defined by
\begin{align}\label{eq-dualzushg}
\hmet{\psi }{\szushg ^*(\psi ')}_0 &:=
\dif \hmet{\psi }{\psi '}_0-\hmet{\szushg \psi }{\psi '}_0,
\qquad \psi ,\psi '\in \mc{S},
\end{align}
is a left connection on $\mc{S}$. It is called the \textit{connection dual to}
$\szushg $ (see also \cite{a-BDMS1}).
The connection dual to $\szushg ^*$, denoted by $\szushg ^{**}$, is again
$\szushg $ itself. Indeed, applying the involution $*$ onto
(\ref{eq-dualzushg}) and having Remark 1 after Definition \ref{d-hmet}
in mind we obtain
\begin{align}
\hmet{\szushg ^*(\psi ')}{\psi }_0 &=
\dif \hmet{\psi '}{\psi }_0-\hmet{\psi '}{\szushg \psi }_0,
\qquad \psi ,\psi '\in \mc{S}.
\end{align}
Hence $\hmet{\psi '}{\szushg \psi }_0=\hmet{\psi '}{\szushg ^{**}\psi }_0$
for all $\psi ,\psi '\in \mc{S}$.

\begin{bem}
It is not clear whether the connection dual to $\szushg $ is linear
or/and compatible with the left multiplication of the quantum Clifford
algebra. The examples in Section \ref{sec-example} show that this can be the
case, but $\szushg $ and $\szushg ^*$ are not necessarily compatible with the
same linear connection on $\Cl (\Gam ,\sigma ,g)$. 
\end{bem}

The mapping $\nabla ^*\szushg :\mc{S}\to \mc{S}$, given by
\begin{align}\label{eq-connlapl}
\hmet{\nabla ^*\szushg (\psi )}{\psi '} &:=
-h g(\hmet{\szushg \psi }{\szushg \psi '}_0),
\qquad \psi ,\psi '\in \mc{S},
\end{align}
is called the \textit{connection Laplacian} associated to the connection
$\szushg $ on $\mc{S}$.

\section{Dirac operators}
\label{sec-Dirac}

Let $\mc{C}$ be a finite dimensional left-covariant $\A $-bimodule with
basis $\{\cb _i\,|\,i=1,\ldots ,p\}$ of $\linv{\mc{C}}$.
Let $\mc{X}=\Lin\{\x _j\,|\,j=1,\ldots ,n\}$ and $\Xalg $
denote the tangent space of a left-covariant FODC $(\Gam ,\dif )$ over $\A $
and the unital complex subalgebra of $\Anull $ generated by $\mc{X}$,
respectively. Observe that $\Xalg $ can be equipped with a filtration
such that $\deg \x _j=1$ for all $j=1,\ldots ,n$.

\begin{defin}\label{d-invdo}
A mapping $\partial :\mc{C}\to \mc{C}$ is called a \textit{left-invariant
differential operator} on $\mc{C}$ (with respect to the differential calculus
$(\Gam ,\dif )$) if there exist functionals
$p_{i,j}\in \Xalg $, $i,j=1,\ldots ,p$, such that
$\partial (a\cb _i)=a_{(1)}p_{j,i}(a_{(2)})\cb _j$
for $i=1,\ldots ,p$. A left-invariant differential operator $\partial $ on
$\mc{C}$ is called an $m^\mathrm{th}$ \textit{order differential operator}
if $\mathrm{deg}\,p_{i,j}\leq m$ for all $i,j=1,\ldots ,p$
with respect to the given filtration of $\Xalg $.
\end{defin}

Let $(\Gam ,\dif )$ be a FODC over $\A $,
$\Cl (\Gam ,\sigma ,g)$ a corresponding quantum Clifford algebra and
$\mc{S}$ a spinor module. Let $\nabla $ and $\nabla '$ be
linear left connections on $\Cl (\Gam ,\sigma ,g)$ and $\mc{S}$, respectively,
which are compatible with the multiplication of the quantum Clifford algebra.
Then we define the Dirac operator $D$ on $\Cl (\Gam ,\sigma ,g)$ by
$D:=\multcl \nabla $ and on $\mc{S}$ by $D:=\multcls \nabla '$.
Equation (\ref{eq-connection}) gives that the Dirac operators on
$\Cl (\Gam ,\sigma ,g)$ and on $\mc{S}$ are first order left-invariant
differential operators. At the end of this section it will be shown that
the connection Laplacian on $\mc{S}$ is a second order differential operator.

\begin{lemma}\label{l-hx}
For all $a,b\in \A $ and $\x \in \mc{X}$ the equation
\begin{align}
h(a_{(1)}\x (a_{(2)})b^*)&=
h\bigl(a(b_{(1)}S^2(\x )^*(b_{(2)}))^*\bigr)
\end{align}
holds.
\end{lemma}

\begin{bew}
Since both sides of the equation are linear in $\x $ it suffices to prove
the lemma for $\x =\x _i$, $i=1,\ldots ,\dim \mc{X}$. We obtain
\begin{align*}
h(a_{(1)}b^*)\x _i(a_{(2)})&=h(a_{(1)}b^*_{(1)})\x _k(a_{(2)})
f^k_l(b^*_{(2)})S(f^l_i)(b^*_{(3)})\\
&=h(a_{(1)}b^*_{(1)})(\x _l(a_{(2)}b^*_{(2)})-\koun (a_{(2)})\x _l(b^*_{(2)}))
S(f^l_i)(b^*_{(3)})\\
&=(h\x _l)(ab^*_{(1)}S(f^l_i)(b^*_{(2)}))
+h(ab^*_{(1)})S(\x _i)(b^*_{(2)}).
\end{align*}
Since the Haar functional is right-invariant, we have
$h\x _l(c)=h(c)\x _l(1)=0$ for $c\in \A $ and therefore
the first summand of the last expression vanishes. On the other hand,
formula (\ref{eq-stern}) gives
\begin{align*}
S(\x _i)(b^*)&=\overline{S^2(\x _i)^*(b)}
\end{align*}
from which the assertion follows.
\end{bew}

Suppose that $\hmet{\cdot }{\cdot }:\mc{S}\ot \cconj{\mc{S}}\to \comp $
is a hermitean metric on $\mc{S}$. Recall that $\beta $ is an involution
of the algebra $\Cl (\Gam ,\sigma ,g)$.

\begin{satz}\label{s-symmdirac}
The Dirac operator $D$ on $\mc{S}$ is symmetric with respect to the metric
$\hmet{\cdot }{\cdot }$ if and only if
$\hmet{D\varphi }{\psi }=\hmet{\varphi }{D\psi }$ for any
$\varphi ,\psi \in \linv{\mc{S}}$ and
equation (\ref{eq-hmetbeta}) is fulfilled
for $u\in \Cl (\Gam ,\sigma ,g)$ and $\varphi ,\psi \in \mc{S}$.
\end{satz}

\begin{bew}
Let $a,b\in \A $ and $\varphi ,\psi \in \linv{\mc{S}}$. Then
\begin{align*}
&\hmet{D(a\varphi )}{b\psi }-\hmet{a\varphi }{D(b\psi )}=\\
&\quad =\hmet{a_{(1)}\x _i(a_{(2)})\w _i\varphi +a\,D\varphi }{b\psi }
-\hmet{a\varphi }{b_{(1)}\x _j(b_{(2)})\w _j\psi +b\,D\psi }\\
&\quad =h(a_{(1)}b^*)\x _i(a_{(2)})\hmet{\w _i\varphi }{\psi }
-h(ab^*_{(1)})\overline{\x _j(b_{(2)})}\hmet{\varphi }{\w _j\psi }
\tag{$*$}\\
&\quad \phantom{=}+h(ab^*)(\hmet{D\varphi }{\psi }-\hmet{\varphi }{D\psi })
\end{align*}
by Definition \ref{d-hmet} and formulas (\ref{eq-connection}) and
(\ref{eq-diffmap}).
Applying Lemma \ref{l-hx} and using formula (\ref{eq-stern}) it follows
that ($*$) is equal to
\begin{align*}
&h(ab^*_{(1)})\overline{B^i_j}\,\overline{\x _j(b_{(2)})}
\hmet{\w _i\varphi }{\psi }
-h(ab^*_{(1)})\overline{\x _j(b_{(2)})}\hmet{\varphi }{\w _j\psi }\\
&\quad +h(ab^*)(\hmet{D\varphi }{\psi }-\hmet{\varphi }{D\psi })\\
&=h(ab^*_{(1)})\overline{B^i_j}\,\overline{\x _j(b_{(2)})}
(\hmet{\w _i\varphi }{\psi }
-\hmet{\varphi }{\beta (\w _i)\psi })
+h(ab^*)(\hmet{D\varphi }{\psi }-\hmet{\varphi }{D\psi }).\tag{$**$}
\end{align*}

Suppose that
$\hmet{D\varphi }{\psi }=\hmet{\varphi }{D\psi }$ for any
$\varphi ,\psi \in \linv{\mc{S}}$ and (\ref{eq-hmetbeta}) is fulfilled.
Since any element of $\mc{S}$ is a linear combination of elements
$a\varphi $, $a\in \A $, $\varphi \in \linv{\mc{S}}$, we obtain
$\hmet{D\varphi }{\psi }=\hmet{\varphi }{D\psi }$ for any
$\varphi ,\psi \in \mc{S}$ by ($**$).

For the other direction of the assertion of the lemma we suppose that
$\hmet{D(a\varphi )}{b\psi }-\hmet{a\varphi }{D(b\psi )}=0$
for $a,b\in \A $ and $\varphi ,\psi \in \linv{\mc{S}}$. Then
trivially $\hmet{D\varphi }{\psi }=\hmet{\varphi }{D\psi }$ for any
$\varphi ,\psi \in \linv{\mc{S}}$. Moreover, from $(**)=0$ we obtain
\begin{align*}
h\bigl(ab^*_{(1)}\overline{B^i_j}\,\overline{\x _j(b_{(2)})}
(\hmet{\w _i\varphi }{\psi }
-\hmet{\varphi }{\beta (\w _i)\psi })\bigr)=0
\end{align*}
for $a,b\in \A $ and $\varphi ,\psi \in \linv{\mc{S}}$.
Since the Haar functional ist left-regular,
we conclude that 
$(b_{(1)}B^i_j\x _j(b_{(2)}))^*(\hmet{\w _i\varphi }{\psi }
-\hmet{\varphi }{\beta (\w _i)\psi })=0$
for any $b\in \A $ and $\varphi ,\psi \in \linv{\mc{S}}$.
Evaluating $\koun $ on this expression and setting $b:=b_k$,
where $\x _i(b_k)=\delta ^i_k$,
we obtain
$\hmet{\w _k\varphi }{\psi }-\hmet{\varphi }{\beta (\w _k)\psi }=0$
for any $\varphi ,\psi \in \linv{\mc{S}}$. Since the elements of $\A $
and the set $\{\w _k\}$ generate $\Cl (\Gam ,\sigma ,g)$,
(\ref{eq-hmetbeta}) holds for any $u\in \Cl (\Gam ,\sigma ,g)$ and
$\varphi ,\psi \in \linv{\mc{S}}$. Finally, by Remark 4 after Definition
\ref{d-hmet} we have
\begin{align*}
\hmet{u(a\varphi )}{b\psi }=\hmet{(\beta (b)ua)\varphi }{\psi }
=\hmet{\varphi }{\beta (a)\beta (u)b\psi }=\hmet{a\varphi }{\beta (u)b\psi }
\end{align*}
for any $a,b\in \A $, $\varphi ,\psi \in \linv{\mc{S}}$ and
$u\in \Cl (\Gam ,\sigma ,g)$.
\end{bew}

\begin{satz}\label{s-laplso}
The connection Laplacian $\nabla ^*\szushg $ associated to the connection
$\szushg $ on $\mc{S}$ is a second order left-invariant differential operator
on $\mc{S}$.
\end{satz}

\begin{bew}
By the defining equation (\ref{eq-connlapl}),
\begin{align}\label{eq-connlaplA0}
\hmet{\nabla ^*\szushg (a\psi )}{b\psi '}&=
-hg\bigl(\hmet{\szushg (a\psi )}{b_{(1)}\x _k(b_{(2)})\w _k\otA \psi '
+b\szushg (\psi ')}_0\bigr)
\end{align}
for $\psi ,\psi '\in \linv{\mc{S}}$, $a,b\in \A $.
Let us reformulate both summands of the expression on the right hand side.
Using (\ref{eq-gSquad}) and Lemma \ref{l-hx} we obtain
\begin{align*}
&hg(\hmet{a\w _i\otA \psi _j}{b_{(1)}\x _k(b_{(2)})\w _k\otA \psi '}_0)=
h(a(b_{(1)}\x _k(b_{(2)}))^*)g(\w _i,\w _k^*)\hmet{\psi _j}{\psi '}\\
&\quad =h(a(b_{(1)}\x _k(b_{(2)}))^*)g(-\overline{B^l_k}\w _l,\w _i)
\hmet{\psi _j}{\psi '}\\
&\quad =-h(a(b_{(1)}B^l_k\x _k(b_{(2)}))^*)g(\w _l,\w _i)
\hmet{\psi _j}{\psi '}\\
&\quad =-h(a_{(1)}\x _l(a_{(2)})b^*)g(\w _l,\w _i)\hmet{\psi _j}{\psi '}
=-\hmet{a_{(1)}\x _l(a_{(2)})g(\w _l,\w _i)\psi _j}{b\psi '}.
\end{align*}
Further, applying (\ref{eq-dualzushg}) we get
\begin{align*}
&hg(\hmet{a\w _i\otA \psi _j}{b\szushg (\psi ')}_0)
=hg(a\w _i \otA \hmet{\psi _j}{\szushg (\psi ')}_0b^*)\\
&\quad =-hg(a\w _i \otA \hmet{\szushg ^*\psi _j}{\psi '}_0b^*)
=-hg(\hmet{a\w _i\otA \szushg ^*\psi _j}{b\psi '}_0).
\end{align*}
{}From this, equation (\ref{eq-connlaplA0}) and since $\hmet{\cdot }{\cdot }$
is non-degenerate we conclude that
\begin{align}\notag
\nabla ^*\szushg (a\psi )=&a_{(1)}\x _k\x _l(a_{(2)})g(\w _k,\w _l)\psi
+a(g\otA \id )(\id \otA \szushg ^*)\szushg (\psi )\\
\label{eq-connlaplexp}
&+a_{(1)}\x _k(a_{(2)})(g\otA \id )(\w _k\otA (\szushg +\szushg ^*)(\psi ))
\end{align}
for all $a\in \A $ and $\psi \in \linv{\mc{S}}$.
\end{bew}

In classical differential geometry the connection Laplacian and the Dirac
operator (corresponding to the torsion-free euclidean connection)
are related by the famous theorem of Bochner \cite{b-LawMich89}. In the
present setting, without any requirement on the torsion,
only a weak form of this assertion holds. For a stronger result
for the quantum group $\SLq 2$ see also Theorem \ref{t-SLq2Bochner}.

\begin{satz}
Let $\szushg $ be a linear left connection on $\mc{S}$ and let
$D$ and $\nabla ^*\szushg $ be the corresponding Dirac operator and
connection Laplacian, respectively. Assume that
$\ker (\id -\sigma )=\ker (\id -\sigma )^2(\subset \Gam \otA \Gam )$. Then
the operator $\nabla ^*\szushg -D^2$ is a first order left-invariant
differential operator on $\mc{S}$.
\end{satz}

\begin{bew}
Since $D=\multcls \szushg $, from (\ref{eq-connection}) it follows that
the operator $D^2$ acts on $\mc{S}$ by the formula
\begin{align*}
D^2(a\psi )&=D(a_{(1)}\x _j(a_{(2)})\w _j\psi +aD\psi )\\
&=a_{(1)}\x _i\x _j(a_{(2)})\w _i\w _j\psi
+a_{(1)}\x _i(a_{(2)})(D(\w _i\psi )+\w _iD\psi )+aD^2\psi
\end{align*}
for $a\in \A ,\psi \in \mc{S}$.
On the other hand, $\nabla ^*\szushg $ satisfies the equation
(\ref{eq-connlaplexp}).
Hence by Definition \ref{d-invdo} we have to show that
the mapping $\partial :\mc{S}\to \mc{S}$, defined by
\begin{align*}
\partial (a\psi )&=a_{(1)}\x _i\x _j(a_{(2)})
(g(\w _i,\w _j)\psi -\w _i\w _j\psi )
\end{align*}
is a first order left-invariant differential operator on $\mc{S}$.

Since $\linv{\Gam }\ot \linv{\Gam }$ is a finite dimensional vector space,
there exists a polynomial $T\in \comp [s]$ in one variable $s$ such that
$(\id -\sigma )T(\sigma )=0$. {}From the assumption 
$\ker (\id -\sigma )=\ker (\id -\sigma )^2$ we conclude that $T(1)\not= 0$.
Hence there exists a polynomial $T'\in \comp [s]$ such that
$1=T(s)/T(1)+(1-s)T'(s)$. Now we have
\begin{align*}\tag{$*$}
\x _i\x _j(a)\w _i\w _j\psi &=\x _i\x _j(a)\bigl(T(\sigma )/T(1)
+(\id -\sigma )T'(\sigma )\bigr)_{ij}^{kl}\w _k\w _l\psi .
\end{align*}
Definition \ref{d-qclifalg} gives that
$T(\sigma )_{ij}^{kl}\w _k\w _l=T(\sigma )_{ij}^{kl}g(\w _k,\w _l)
=T(1)g(\w _i,\w _j)$ because of $g\sigma =g$. On the other hand,
the second requirement on the braiding on page \pageref{breq2}
ensures that $\x _i\x _j(\id -\sigma )_{ij}^{kl}=[\x _k,\x _l]$
are elements of $\mc{X}$ for all $k,l$. Therefore, equation ($*$) gives that
\begin{align*}
\x _i\x _j(a)\w _i\w _j\psi &=\x _i\x _j(a)g(\w _i,\w _j)\psi
+[\x _m,\x _n]T'(\sigma )_{mn}^{kl}\w _k\w _l\psi .
\end{align*}
This means that $\partial $ is a first order left-invariant
differential operator on $\mc{S}$.
\end{bew}

\section{The quantum group {\boldmath $\SLq 2$}. An example.}
\label{sec-example}

Let $\A $ denote the Hopf algebra $\OSLq 2$ with generators $u^i_j$.
We use the symbols $\Rda =(\Rda ^{ij}_{kl})$, $\oC =(\oC ^{ij})$ and
$\uC =(\uC _{ij})$ for the matrices with entries
\begin{equation}\label{eq-RCmatrix}
\begin{aligned}
\Rda ^{ij}_{kl}&=\delta ^i_l\delta ^j_k q^{\delta ^i_j-1/2}
+q^{-1/2}\qm \delta ^i_1\delta ^k_1\delta ^j_2\delta ^l_2,\\
\oC ^{ij}&=q^{-1/2}\delta ^i_1\delta ^j_2-q^{1/2}\delta ^i_2\delta ^j_1,\
\uC _{ij}=-\oC ^{ij},
\end{aligned}
\end{equation}
$i,j,k,l=1,2$. Further notations are collected in Appendix \ref{sec-appslq2}.

Let $(\Gam ,\dif )$ be one of the
($S^2$-invariant bicovariant) $4D_\pm $-calculi over $\A $.
There is a basis
$\{\w _{ij}\,|\,i,j=1,2\}$ of $\linv{\Gam }$ such that
$\rkow (\w _{ij})=\w _{kl}\ot u^k_iu^l_j$ and
$\w _{ij}\ract a:=S(a_{(1)})\w _{ij}a_{(2)}=f^{ij}_{kl}(a)\w _{kl}$, where
\begin{align}\label{eq-fmatrix}
f^{ij}_{kl}&=\koun _\pm \begin{pmatrix}
1 & 0 & q^{-1/2}\qm FK & 0\\
-q^{1/2}\qm K^{-1}E & K^{-2} & -\qm ^2FE & q^{-1/2}\qm FK^{-1}\\
0 & 0 & K^2 & 0\\
0 & 0 & -q^{1/2}\qm KE & 1
\end{pmatrix}.
\end{align}
Equivalently, the right multiplication of the differential 1-forms $\w _{ij}$
by the generators $u^k_l$ of $\A $ are given by
\begin{align}\label{eq-wuvert}
\w _{ij}u^k_l&=\spinsig{0}\Rda {}^{mr}_{it}\Rdam {}^{ts}_{jl}u^k_m\w _{rs},
\end{align}
where $\spinsig{0}=\pm $, depending on the sign of the differential calculus.
The element $\w :=q^{1/2}\oC^{ij}/\qm \,\w _{ij}$ of $\Gam $ is biinvariant
and the differential $\dif :\A \to \Gam $
can be defined by $\dif a:=\w a-a\w $ for any $a\in \A $.
The dual basis to $\{\w _{ij}\}$ of the quantum tangent space $\mc{X}$ is
\begin{equation}
\begin{aligned}
\{
&\x _{11}=-q^{1/2}\koun _\pm K^{-1}E,
\x _{12}=1/\qm \,(\koun _\pm K^{-2}-\koun ),\\
&\x _{21}=-\qm \koun _\pm FE -q/\qm (\koun _\pm K^2-\koun ),
\x _{22}=q^{-1/2}\koun _\pm FK^{-1}
\}
\end{aligned}
\end{equation}
For these functionals we have
\begin{align}
S^2(\x _{11})&=q^2\x _{11},&
S^2(\x _{12})&=\x _{12},&
S^2(\x _{21})&=\x _{21},&
S^2(\x _{22})&=q^{-2}\x _{22}.
\end{align}

We consider three involutions $\ast $
of $\A $. Let $\invcsu $ denote the involution of $\OSUq 2$,
$\invncsu $ the one of
$\OSUq{1,1}$ and $\invslR $ the one of $\OSLq{2,\mathbb{R}}$. In all three
cases the involution of $\A $ can be uniquely extended to an involution
$\ast $ of $\Gam $ such that $\dif (a^\ast )=(\dif a)^\ast $ for any
$a\in \A $. Then the involution of $\mc{X}$ is given by
\begin{equation}
\begin{aligned}
\x _{11}^\invcsu &=-q\x _{22},\;& \x _{12}^\invcsu &=\x _{12},\;&
 \x _{21}^\invcsu &=\x _{21},\;& \x _{22}^\invcsu &=-q^{-1}\x _{11},\\
\x _{11}^\invncsu &=q\x _{22},\;& \x _{12}^\invncsu &=\x _{12},\;&
 \x _{21}^\invncsu &=\x _{21},\;& \x _{22}^\invncsu &=q^{-1}\x _{11},\\
\x _{11}^\invslR &=\x _{11},\;& \x _{12}^\invslR &=-\x _{12},\;&
 \x _{21}^\invslR &=-\x _{21}{-}\qm \x _{12},\;&
 \x _{22}^\invslR &=q^2\x _{22}.
\end{aligned}
\end{equation}
Hence the matrices $B=(B^{ij}_{kl})$, defined by
$S^2(\x _{ij})^*=B^{ij}_{kl}\x _{kl}$, take for the involutions
$\invcsu $, $\invncsu $ and $\invslR $ the form
\begin{align}\label{eq-bmat}
\begin{pmatrix}
0 & 0 & 0 & -q^3\\
0 & 1 & 0 & 0\\
0 & 0 & 1 & 0\\
-q^{-3} & 0 & 0 & 0
\end{pmatrix}, \quad
\begin{pmatrix}
0 & 0 & 0 & q^3\\
0 & 1 & 0 & 0\\
0 & 0 & 1 & 0\\
q^{-3} & 0 & 0 & 0
\end{pmatrix} \text{ and }
\begin{pmatrix}
q^{-2} & 0 & 0 & 0\\
0 & -1 & 0 & 0\\
0 & -\qm & -1 & 0\\
0 & 0 & 0 & q^4
\end{pmatrix},
\end{align}
respectively.
For the involution of $\Gam $ we obtain the formulas
\begin{equation}\label{eq-winv}
\begin{aligned}
\w _{11}^\invcsu &=q^{-1}\w _{22},\;& \w _{12}^\invcsu &=-\w _{12},\;&
 \w _{21}^\invcsu &=-\w _{21},\;& \w _{22}^\invcsu &=q\w _{11},\\
\w _{11}^\invncsu &=-q^{-1}\w _{22},\;& \w _{12}^\invncsu &=-\w _{12},\;&
 \w _{21}^\invncsu &=-\w _{21},\;& \w _{22}^\invncsu &=-q\w _{11},\\
\w _{11}^\invslR &=-\w _{11},\;& \w _{12}^\invslR &=\w _{12}{-}\qm \w _{21},\;&
 \w _{21}^\invslR &=\w _{21},\;& \w _{22}^\invslR &=-q^{-2}\w _{22}.
\end{aligned}
\end{equation}

Recall that the braiding $\sigma $ of $\Gam $
(defined by (\ref{eq-bikovbraid})) is given by the matrix
\begin{align}\label{eq-braiding}
\sigma ^{ij,kl}_{mn,rs}&=f^{mn}_{kl}(u^i_ru^j_s)=(\Rda _{23}\Rda _{12}
\Rdam _{34}\Rdam _{23})^{ijkl}_{mnrs}\qquad \text{(leg numbering)}
\end{align}
and the setting
$\sigma (\w _{mn}\otA \w _{rs})=\sigma ^{ij,kl}_{mn,rs}\w _{ij}\otA \w _{kl}$.

\begin{lemma}\label{l-bicsigmet}
A linear mapping $g:\Gam \otA \Gam \to \A $ is a bicovariant $\sigma $-metric
of $\Gam $ if and only if it is of the form
\begin{align}\label{eq-bisigmet}
g(\w _{ij}\otA \w _{kl})=-q^{1/2}c_1 \uC _{im}\uC _{nl}\Rdam {}^{mn}_{jk},
\qquad c_1\in \realx .
\end{align}
\end{lemma}

\begin{bew}
Usual methods (see e.\,g.~\cite{a-HeckSchm1}) give that $g$ is a bicovariant
mapping if and only if $g(\w _{ij}\otA \w _{kl})=\lambda \uC _{ij}\uC _{kl}
+\mu \uC _{jk}\uC _{il}$, $\lambda ,\mu \in \comp $. Equation
(\ref{eq-wuvert}) implies that $g$ is of the form (\ref{eq-bisigmet}) with
$c_1\in \comp $.
Using (\ref{eq-winv}), from the compatibility with the involution we conclude
in all three cases that $c_1\in \real $.
Finally, non-degeneracy of $g$ gives that $c_1\in \realx $.
The proof of the property $g\sigma =g$ is an easy computation. The other
requirements are fulfilled for each homomorphism of the the bicovariant
$\A $-bimodules $\Gam \otA \Gam $ and $\A $ (see \cite{a-HeckSchu00}).
\end{bew}

Let us now fix such a $\sigma $-metric $g$. The nonzero entries of the matrix
$(g_{ij,kl})$, $g_{ij,kl}:=g(\w _{ij}\otA \w _{kl})$ are
\begin{align}
g_{11,22}&=-c_1,&
g_{12,21}&=c_1,&
g_{21,12}&=c_1,&
g_{22,11}&=-q^2c_1,&
g_{12,12}&=\qm c_1.
\end{align}
By Definition \ref{d-qclifalg}, the quantum Clifford algebra
$\Cl (\Gam ,\sigma ,g)$ is generated by the following set of relations:
\begin{equation}\label{eq-cliffrel}
\begin{aligned}
\{&
 \w _{11}^2, 
 \w _{12}\w _{11}+q^2\w _{11}\w _{12}+\qm \w _{11}\w _{21},
 \w _{21}\w _{11}+\w _{11}\w _{21},
 \w _{22}\w _{11}+\w _{11}\w _{22}+(q^2+1)c_1,\\
& \w _{12}^2+\qm \w _{11}\w _{22},
 \w _{21}\w _{12}+\w _{12}\w _{21}+q^{-1}\qm \w _{11}\w _{22}-(1+q^{-2})c_1,\\
& \w _{22}\w _{12}+q^2\w _{12}\w _{22}+\qm \w _{21}\w _{22},
 \w _{21}^2,
 \w _{22}\w _{21}+\w _{21}\w _{22},
 \w _{22}^2
\}.
\end{aligned}
\end{equation}
Setting $|\w _{ij}|=3-i-j$ and $|u^k_l|=k-l$, from equations (\ref{eq-wuvert})
and (\ref{eq-cliffrel}) we directly obtain that
the algebra $\Cl (\Gam ,\sigma ,g)$ becomes a graded algebra with grading
$|\cdot |$.

\begin{satz}\label{s-clsimple}
The quantum Clifford algebra $\Cl (\Gam ,\sigma ,g)$ is a simple
left-covariant algebra.
The elements
\begin{equation}
\begin{aligned}
\psi ^+_1&:=\frac{q}{\qp c_1}\w _{11}\w _{12}\w _{21}\w _{22},\quad &
\psi ^+_2&:=\w _{21}\w _{22},\\
\psi ^-_1&:=\w _{11}\w _{21}\w _{22},&
\psi ^-_2&:=q\w _{12}\w _{21}\w _{22}
\end{aligned}
\end{equation}
generate a minimal left-covariant left ideal $\mc{S}_0$ of
$\Cl (\Gam ,\sigma ,g)$.
\end{satz}

The left-covariant left ideals $\A \psi ^\eta _1 \oplus \A \psi ^\eta _2$,
$\eta \in \{+,-\}$, are denoted by $\mc{S}^\eta _0$. 

\begin{bew}[ of the Proposition]
For the first assertion it is enough to show that there are no non-trivial
left-covariant ideals $\mc{J}$ of $\Cl (\Gam ,\sigma ,g)$.

Suppose that $\mc{J}$ is a nonzero left-covariant ideal of
$\Cl (\Gam ,\sigma ,g)$. Then $\mc{J}$ is a free left $\A $-module with
a left-invariant basis. Since $\Cl (\Gam ,\sigma ,g)$ is finite dimensional,
$\mc{J}$ is so, too.
Let $\rho $ be a nonzero left-invariant element of $\mc{J}$. We show that then
$1\in \mc{J}$, from which $\mc{J}=\Cl (\Gam ,\sigma ,g)$ follows.
\begin{enumerate}
\item
If $\w _{11}\rho =0$ then $\rho =\w _{11}\rho '$ for a
$\rho '\in \Cl (\Gam ,\sigma ,g)$. Otherwise $\w _{11}\rho \not= 0$.
Hence there exists $\rho _1\in \mc{J}\setminus \{0\}$
such that $\rho _1=\w _{11}\rho '$.
\item
If $\rho _1\w _{11}=0$ then $|\rho _1|=1$.
Otherwise $\rho _1\w _{11} \not= 0$ and we get $|\rho _1\w _{11}|=1$.
Hence there exists $\rho _2\in \mc{J}\setminus \{0\}$
such that $|\rho _2|=1$.
\item
Let $\rho _2\in \mc{J}\setminus \{0\}$, $|\rho _2|=1$.
Then $\rho _2$ is a linear combination of the elements $\w _{11}$,
$\w _{11}\w _{12}$, $\w _{11}\w _{21}$ and $\w _{11}\w _{12}\w _{21}$.
If $\rho _2\theta _{21}=0$ then there are $\lambda _1,\lambda _2\in \comp $
such that $\rho _3:=\rho _2=\w _{11}(\lambda _1+\lambda _2\w _{12})\w _{21}$.
Otherwise $\rho _2\w_{21}\not= 0$ and there are
$\lambda _1,\lambda _2\in \comp $ such that
$\rho _3:=\rho _2\w _{21}=\w _{11}(\lambda _1+\lambda _2\w _{12})\w _{21}$.
Again $\rho _3\in \mc{J}\setminus \{0\}$.
\item
We have $\w _{21}\rho _3=-\lambda _2(1+q^{-2})c_1\w _{11}\w _{21}$.
If $\w _{21}\rho _3=0$ then $\lambda _2=0$ and hence $\rho _3$ is a nonzero
multiple of $\w _{11}\w _{21}$.
Otherwise $\w _{21}\rho _3\not=0$ and therefore
$\w _{21}\rho _3$ is a nonzero multiple of $\w _{11}\w _{21}$.
\item
We obtained that $\rho _4:=\w _{11}\w _{21}\in \mc{J}$.
Then $\mc{J}\ni \w _{12}\rho _4=-q^2\w _{11}\w _{12}\w _{21}$.
We also have $\mc{J}\ni \rho _4\w _{12}=-\w _{11}\w _{12}\w _{21}
+(1+q^{-2})c_1\w _{11}$, hence $\w _{11}\in \mc{J}$.
Therefore $\w _{11}\w _{22}\in \mc{J}$ and $\mc{J}\ni \w _{22}\w _{11}
=-\w _{11}\w _{22}-(1+q^2)c_1$. Hence $1\in \mc{J}$.
\end{enumerate}

Now we turn to the proof of the second assertion.
Using (\ref{eq-cliffrel}) one can easily see
that $\Cl (\Gam ,\sigma ,g)\mc{S}_0\subset \mc{S}_0$.
Moreover, since the dimension of the algebra of left-invariant
elements of $\Cl (\Gam ,\sigma ,g)$ is 16, each minimal left-covariant
left ideal of $\Cl (\Gam ,\sigma ,g)$ is 4-dimensional.
\end{bew}

With help of the relations (\ref{eq-cliffrel}) of $\Cl (\Gam ,\sigma ,g)$
one can easily determine the left action of the generators $\w _{ij}$
on $\mc{S}_0$. We obtain the following table.
\begin{displaymath}
\begin{array}{|r||l|l|l|l|}
\hline
& \psi ^+_1 & \psi ^+_2 & \psi ^-_1 & \psi ^-_2 \\
\hline
\w _{11} & 0 & \psi ^-_1 & 0 & \qp c_1\psi ^+_1 \\
\w _{12} & -\qm \psi ^-_1 & q^{-1}\psi ^-_2 & -q\qp c_1\psi ^+_1 & 0\\
\w _{21} & -\psi ^-_1 & 0 & 0 & \qp c_1\psi ^+_2 \\
\w _{22} & -q\psi ^-_2 & 0 & -q\qp c_1\psi ^+_2 & 0\\
\hline
\end{array}
\end{displaymath}
Equivalently we can write
\begin{align}\label{eq-clwirk}
\w _{ij}\psi ^+_k&=-\Rda {}^{lm}_{ij}\uC _{mk}\psi ^-_l, &
\w _{ij}\psi ^-_k&=-q^{1/2}\qp c_1\psi ^+_i\uC _{jk}.
\end{align}

The following lemma proves that $\mc{S}_0$ is a spinor module of
$\Cl (\Gam ,\sigma ,g)$.

\begin{lemma}\label{l-spinbim}
The left $\A $-modules $\mc{S}_0,\mc{S}^+_0$ and $\mc{S}^-_0$ are invariant
under right multiplication by $\A $.
\end{lemma}

\begin{bew}
By Proposition \ref{s-clsimple}, $\mc{S}_0=\Cl (\Gam ,\sigma ,g)\psi ^+_2$.
Since
\begin{align*}
\psi ^+_2 a=& \w _{21}\w _{22}a
=a_{(1)}f^{21}_{ij}f^{22}_{kl}(a_{(2)})\w _{ij}\w _{kl}\\
=&-q^{1/2}\qm a_{(1)}K^2KE(a_{(2)})\w _{21}\w _{21}
 +a_{(1)}K^2(a_{(2)})\w _{21}\w _{22}
=a_{(1)}K^2(a_{(2)})\psi ^+_2
\end{align*}
for all $a\in \A $, we have $(u\psi ^+_2)a=ua_{(1)}K^2(a_{(2)})\psi ^+_2
\in \Cl (\Gam ,\sigma ,g)\psi ^+_2=\mc{S}_0$ for all
$u\in \Cl (\Gam ,\sigma ,g)$ and $a\in \A $.
Further, $\mc{S}^\eta _0=\Cl ^\eta (\Gam ,\sigma ,g)\psi ^+_2$
for $\eta \in \{+,-\}$. This implies
that $(u\psi ^+_2)a=ua_{(1)}K^2(a_{(2)})\psi ^+_2\in
\Cl ^\eta (\Gam ,\sigma ,g)\psi ^+_2=\mc{S}^\eta _0$ for all
$u\in \Cl ^\eta (\Gam ,\sigma ,g)$ and $a\in \A $.
\end{bew}

By Lemma \ref{l-spinbim} there exists a matrix
$\tilde{f}=(\tilde{f}^i_j)_{i,j=1,\ldots ,4}$, $\tilde{f}^i_j\in \Anull $,
such that
$\psi ^+_i a=a_{(1)}(\tilde{f}^i_j(a_{(2)})\psi ^+_j
+\tilde{f}^i_{j+2}(a_{(2)})\psi ^-_j)$ and
$\psi ^-_i a=a_{(1)}(\tilde{f}^{i+2}_j(a_{(2)})\psi ^+_j
+\tilde{f}^{i+2}_{j+2}(a_{(2)})\psi ^-_j)$
for any $a\in \A $ and $i=1,2$, where the sum is running over $j=1,2$.
Similarly to the proof of the lemma one can compute the matrix elements
$\tilde{f}^i_j$.
Using the commutation relations (\ref{eq-cliffrel}) of $\Cl (\Gam ,\sigma ,g)$
and (\ref{eq-fmatrix}) one gets
\begin{align}\label{eq-spinfmatrix}
\tilde{f}^i_j&=\begin{pmatrix}
1 & q^{-1/2}\qm FK & 0 & 0\\
0 & K^2 & 0 & 0\\
0 & 0 & \koun _\pm K^2 & 0\\
0 & 0 & -q^{-1/2}\qm \koun _\pm KE & \koun _\pm
\end{pmatrix}.
\end{align}
For the generators $u^i_j$ of $\A $ this means that
\begin{align}\label{eq-psiuvert}
\psi ^+_iu^j_k&=\spinsig{+}q^{k-3/2}u^j_l\psi ^+_m\Rda {}^{lm}_{ik},&
\psi ^-_iu^j_k&=\spinsig{-}q^{k-3/2}u^j_l\psi ^-_m\Rdam {}^{lm}_{ik},
\end{align}
where the notation $\spinsig{+}:=1$, $\spinsig{-}:=\spinsig{0}=\pm 1$
(depending on the sign of the differential calculus) is used.

Let $\chi $ denote the algebra automorphism of $\A $ defined by
$\chi \lact u^i_j :=q^{j-i}u^i_j$ ($\chi \lact a=q^{-|a|}a$ for homogeneous
elements $a\in \A $ with respect to the grading $|\cdot |$)
and let $\tilde{\A }$ denote the left crossed product algebra
$\A \rtimes \comp [\chi ,\chi ^{-1}]$
with multiplication $(a\chi ^k)(b\chi ^l)=a(\chi ^k\lact b)\chi ^{k+l}$.
Then obviously $\tilde{\A }$ is an $\A $-bimodule.
The settings $\kopr (\chi ):=\chi \ot \chi $ and $S(\chi )=\chi ^{-1}$
turn $\tilde{\A }$ into a Hopf algebra.
Let us define a right coaction
$\rkow ':\linv{(\mc{S}_0)}\to \linv{(\mc{S}_0)}\ot \tilde{\A }$
on the vector space $\linv{(\mc{S}_0)}$ by
\begin{align}\label{eq-rspinkow}
\rkow '(\psi ^+_i)&=\psi ^+_j\ot u^j_i\chi ,&
\rkow '(\psi ^-_i)&=\psi ^-_j\ot u^j_i\chi ,& i&=1,2.
\end{align}
Recall the definition of a Doi-Hopf module in \cite{a-Doi92}.

\begin{thm}
(i) The right coaction $\rkow '$ on $\linv{(\mc{S}_0)}$ can be extended
to a right coaction $\rkow $ on the left-covariant $\A $-bimodule $\mc{S}_0$.
This means that
$\mc{S}_0$ together with the multiplication from the right
as a right action of $\A $
and $\rkow $ as a right coaction of $\tilde{\A }$ becomes also a Doi-Hopf
module in the category $\mathbf{M} (\A )^{\tilde{\A }}_\A $.\\
(ii) The right coaction $\rkow $ on $\mc{S}_0$ is compatible with the
left multiplication of $\Cl (\Gam ,\sigma ,g)$ on $\mc{S}_0$, i.\,e.\
$\rkow (u\psi )=\rkow (u)\rkow (\psi )$ for any $u\in \Cl (\Gam ,\sigma ,g)$
and $\psi \in \mc{S}_0$.
\end{thm}

This theorem gives reason to think about $\tilde{\A }$ as the function
algebra of the quantum spin group corresponding to the quantum Clifford
algebra $\Cl (\Gam ,\sigma ,g)$.

\begin{bew}
We set $\rkow (a\psi ^\eta _i)=a_{(1)}\psi ^\eta _j\ot a_{(2)}u^j_i\chi $
for $a\in \A $, $\eta \in \{+,-\}$ and $i=1,2$.
Then compatibility of $\rkow $ with the right action of $\A $ on $\mc{S}_0$
means that $\rkow (\psi a)=\rkow (\psi )\kopr (a)$. {}From (\ref{eq-psiuvert})
we conclude that
\begin{align*}
\rkow (\psi ^\eta _iu^j_k)&=\rkow (\spinsig{\eta }q^{k-3/2}u^j_l\psi ^\eta _m
\Rda ^\eta {}^{lm}_{ik})
=\spinsig{\eta }q^{k-3/2}u^j_r\psi ^\eta _s
\Rda ^\eta {}^{lm}_{ik}\ot u^r_lu^s_m\chi \\
&=\spinsig{\eta }q^{k-3/2}u^j_r\psi ^\eta _s
\ot u^l_iu^m_k\Rda ^\eta {}^{rs}_{lm}\chi
=\spinsig{\eta }q^{m-3/2}u^j_r\psi ^\eta _s\Rda ^\eta {}^{rs}_{lm}
\ot q^{k-m}u^l_iu^m_k\chi \\
&=\psi ^\eta _lu^j_m\ot u^l_i\chi u^m_k =\rkow (\psi ^\eta _i)\kopr (u^j_k).
\end{align*}

For the second assertion it suffices to prove the formula
$\rkow (\w _{ij}\psi ^\eta _k)=\rkow (\w _{ij})\rkow (\psi ^\eta _k)$
for $\eta \in \{+,-\}$. Using
(\ref{eq-clwirk}) we obtain
\begin{align*}
\rkow (\w _{ij})\rkow (\psi ^+_k)&=\w _{mn}\psi ^+_l\ot u^m_iu^n_ju^l_k\chi
=-\Rda ^{rs}_{mn}\uC _{sl}\psi ^-_r \ot u^m_iu^n_ju^l_k\chi \\
&=-\uC _{sl}\psi ^-_r\ot u^r_mu^s_n\Rda ^{mn}_{ij}u^l_k\chi
=-\Rda ^{mn}_{ij}\psi ^-_r\ot u^r_m\uC _{nk}\chi \\
&=\rkow (-\Rda ^{mn}_{ij}\uC _{nk}\psi ^-_m)
=\rkow (\w _{ij}\psi ^+_k).
\end{align*}
The statement for $\eta =-$ can be shown similarly.
\end{bew}

One can ask whether the introduced $*$-structures on $\A $ can be extended
to the left-covariant $\A $-bimodule $\mc{S}_0$. Such an extension is not
unique in general. One calls the involutions $*$ and $*'$ of $\mc{S}_0$
\textit{equivalent} if there exists an automorphism $\varphi $ of the
left-covariant $\A $-bimodule $\mc{S}_0$, $\varphi (a\psi b)=a\varphi (\psi )b$
for $a,b\in \A ,\psi \in \mc{S}_0$,
such that $\varphi (\psi ^*)=\varphi (\psi )^{*'}$ for all $\psi \in \mc{S}_0$.
To determine the involutions of $\mc{S}_0$, we use the necessary conditions
\begin{align}
\psi _i^*&=\lambda ^i_j\psi _j,\quad \lambda ^i_j\in \comp ,&
\lambda ^i_j\tilde{f}^j_k&=(\tilde{f}^i_l)^*\lambda ^l_k, \quad i,k=1,2,3,4,
\end{align}
which stem from the compatibility of the involution with the left coaction
and the bimodule structure of $\mc{S}_0$, respectively.
Here we used the notation $\psi _i=\psi _i^+$ for $i=1,2$ and
$\psi _i=\psi _{i-2}^-$ for $i=3,4$.

Let us consider the $4D_+$-calculus. The involutions $(\cdot )^\invcsu $
and $(\cdot )^\invncsu $ of $\A $ can be extended to involutions of the
left-covariant $\A $-bimodule  $\mc{S}_0$. Moreover, these involutions are
unique up to equivalence of left-covariant bimodules. In particular,
the formulas
\begin{equation}
\begin{aligned}
(\psi ^+_1)^\invcsu &=-\psi ^-_2,\,&
(\psi ^+_2)^\invcsu &=\psi ^-_1,\,&
(\psi ^-_1)^\invcsu &=\psi ^+_2,\,&
(\psi ^-_2)^\invcsu &=-\psi ^+_1,\\
(\psi ^+_1)^\invncsu &=\psi ^-_2,\,&
(\psi ^+_2)^\invncsu &=\psi ^-_1,\,&
(\psi ^-_1)^\invncsu &=\psi ^+_2,\,&
(\psi ^-_2)^\invncsu &=\psi ^+_1
\end{aligned}
\end{equation}
hold.
Observe that the setting $\chi ^*:=\chi $ together with each of the three
introduced involutions of $\A $ turn $\tilde{\A }$ into a Hopf $*$-algebra.
It is easy to see that the involutions $\invcsu $ and $\invncsu $
of $\mc{S}_0$ are compatible with the right coaction $\rkow $ of $\tilde{\A }$.
On the other hand, if we consider the $4D_-$-calculus then
there exists no involution of the left-covariant $\A $-bimodule
$\mc{S}_0 $ extending $(\cdot )^\invcsu $ or $(\cdot )^\invncsu $ of $\A $.
However, for both calculi the
involution $(\cdot )^\invslR $ of $\A $ can be extended to an involution of
the left-covariant $\A $-bimodule $\mc{S}_0$.
This involution is unique up to equivalence of left-covariant $\A $-bimodules,
and also compatible with the right coaction of $\tilde{\A }$.
The corresponding formulas read as
\begin{align}
(\psi ^+_1)^\invslR &=q^{-1}\psi ^+_1,&
(\psi ^+_2)^\invslR &=\psi ^+_2,&
(\psi ^-_1)^\invslR &=\psi ^-_1,&
(\psi ^-_2)^\invslR &=q\psi ^-_2.
\end{align}

\begin{thm}\label{t-unmet}
For each one of the given involutions $*$ of $\A $ there exists a
left-covariant (hermitean non-degenerate) metric $\hmet{\cdot }{\cdot }$ on
$\mc{S}_0$ which satisfies (\ref{eq-hmetbeta}). Moreover, this metric is
unique up to a nonzero real factor $c_2$.
The explicit formula for the metric is given by
\begin{equation}\label{eq-unmet}
\begin{aligned}
\hmet{\psi ^\eta _i}{\psi ^{\eta '}_j}&=
c_2\delta ^\eta _{-\eta '}\delta _{i,j}(\delta ^i_1+q^2\delta ^i_2),\\
\hmet{\psi ^\eta _i}{\psi ^{\eta '}_j}&=
c_2\delta ^\eta _{-\eta '}\delta _{i,j}(\delta ^i_1-q^2\delta ^i_2),\\
\hmet{\psi ^\eta _i}{\psi ^{\eta '}_j}&=
c_2\delta ^\eta _{\eta '}\delta _{i+j,3}
(q^{-1}\delta _{i,1}+q\delta _{i,2})
(\delta _{\eta ,+}+\qp c_1\delta _{\eta ,-})
\end{aligned}
\end{equation}
for the involutions $\invcsu $, $\invncsu $ and $\invslR $, respectively.
\end{thm}

Observe that non of the obtained metrics is positive definite.

\begin{bew}
Since $\Cl (\Gam ,\sigma ,g)$ is a simple left-covariant algebra by
Proposition \ref{s-clsimple},
the existence of such a metric follows from Proposition \ref{s-spinmet}.
To obtain the given formulas for the metric, one should choose the elements
$x=-q\qp c_1/c_2\psi _1^-$,
$-q\qp c_1/c_2\psi _1^-$ and
$-q^{-2}\qp c_1/c_2\psi _2^+$ of $\beta (\mc{S})\mc{S}$
in case of the involution
$\invcsu $, $\invncsu $ and $\invslR $, respectively.

Now we prove the uniqueness of the metric for the involution $\invcsu $.
For the other involutions the proof is analogous.
Since $\mc{S}$ is a minimal left
ideal, each vector $\psi \in \linv{\mc{S}}\setminus \{0\}$ is cyclic.
Hence because of (\ref{eq-hmetbeta}) the numbers $\hmet{ \psi _i^\pm }{\psi }$
determine the metric completely. Set $\psi :=\psi ^-_2$.
{}From the table before equation (\ref{eq-clwirk}) we obtain that
$\w _{12}\psi ^-_2=0$. Hence $0=\hmet{\psi '}{\w _{12}\psi ^-_2}=
\hmet{\beta (\w _{12})\psi '}{\psi ^-_2}$. By (\ref{eq-bmat}) it follows
that $\beta (\w _{12})=\w _{12}$. Setting $\psi '=\psi ^+_1$, $\psi ^+_2$
and $\psi ^-_1$, again from the table before (\ref{eq-clwirk}) we conclude
that $\hmet{\psi ^-_1}{\psi ^-_2}=\hmet{\psi ^-_2}{\psi ^-_2}=
\hmet{\psi ^+_1}{\psi ^-_2}=0$. Hermiticity and non-degeneracy of the metric
imply that the remaining parameter $\hmet{\psi ^+_2}{\psi ^-_2}$ is a
nonzero real number.
\end{bew}

Let us define a right coaction $\rkow $ of
$\A \rtimes \comp [\chi ,\chi ^{-1}]$ on $\mc{S}_0\ot \cconj{\mc{S}_0}$ 
by
\begin{align}
\rkow (\psi \ot \cconj{\psi '{}}):=(\psi _{(0)}\ot \cconj{\psi '_{(0)}{}})\ot
\psi _{(1)}\chi ^{-2}(\psi '_{(1)})^*,
\end{align}
where $*$ is the involution of
$\A \rtimes \comp [\chi ,\chi ^{-1}]$. Since $\chi ^*=\chi $,
the image of $\rkow $ is a subspace of $\mc{S}\ot \cconj{\mc{S}}\ot \A $.
It is not difficult to prove that for the involutions $\invcsu $ and
$\invncsu $ (but not for $\invslR $)
the mappings $\hmet{\cdot }{\cdot }_0:\mc{S}_0\ot \cconj{\mc{S}_0}\to \A $ and
$\hmet{\cdot }{\cdot }:\mc{S}_0\ot \cconj{\mc{S}_0}\to \comp $ are
right-covariant:
$(\hmet{\cdot }{\cdot }_0\ot \id )\rkow =\kopr \hmet{\cdot }{\cdot }_0$
on $\mc{S}_0\ot \cconj{\mc{S}_0}$.

\begin{satz}\label{s-lconnslq2}
Let $g$ be one of the $\sigma $-metrics of Lemma \ref{l-bicsigmet}.
Then there exist 8 linear left connections $\nabla $ on $\Gam $
compatible with the multiplication of the quantum Clifford algebra
$\Cl (\Gam ,\sigma ,g)$. They can be given by the formula (\ref{eq-innzushg})
with one of the mappings $\ts :=\nu \ts _{(i)}$, $i=1,2,3,4$, where
$\nu \in \{+,-\}$,
\begin{align}\label{eq-sigmatilde}
\ts _{(i)}(\w _{kl}\ot \w _{mn})&=
(\Rda _{23}\Rda ^{\eta _i}_{12}\Rda ^{\eta '_i}_{34}\Rdam _{23})^{rstx}_{klmn}
\w _{rs}\ot \w _{tx},
\end{align}
and $(\eta _i,\eta '_i)$ denotes the pair of signs $(+,+)$, $(+,-)$, $(-,+)$
and $(-,-)$ for $i=1,2,3$ and $4$, respectively.
\end{satz}

\begin{bem}
None of left connections $\nabla $ in Proposition \ref{s-lconnslq2} satisfies
(\ref{eq-torshom}). For example, we have
\begin{align}
(\id -\sigma )(\id +q\ts _{(1)})&=0, &
(\id -\sigma )(\id +q^{-2}\ts _{(2)})(\id +q^2\ts _{(2)})&=0.
\end{align}
\end{bem}

\begin{bew}[ of the Proposition]
We only have to show that the assumptions of Proposition \ref{s-clconn}
are fulfilled. This is very easily done for $\nu \ts _{(2)}=\nu \sigma $
and $\nu \ts _{(3)}=\nu \sigma ^{-1}$. Let us verify them for
$\nu \ts _{(1)}$. {}From (\ref{eq-wuvert}) follows that $\nu \ts _{(1)}$
defines a homomorphism of the $\A $-bimodule $\Gam \otA \Gam $. Indeed,
the coefficient of $u^r_t\w _{xy}\otA \w _{zw}$ of the expression
$\nu \ts _{(1)}(\w _{kl}\otA \w _{mn})u^r_s-
\nu \ts _{(1)}(\w _{kl}\otA \w _{mn}u^r_s)$ is
\begin{align*}
(\nu \Rda _{12}\Rdam _{23}\Rda _{34}\Rdam _{45}
\Rda _{23}\Rda _{12}\Rda _{34}\Rdam _{23}
-\nu \Rda _{34}\Rda _{23}\Rda _{45}\Rdam _{34}
\Rda _{12}\Rdam _{23}\Rda _{34}\Rdam _{45})^{txyzw}_{klmns},
\end{align*}
or
\begin{gather*}
\epsfig{file=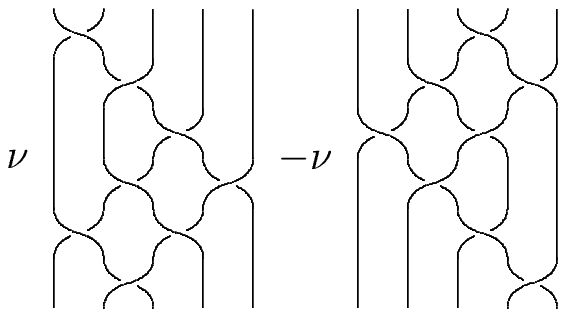}
\end{gather*}
in the graphical calculus. This is zero because of the Yang-Baxter equation.
Inserting the matrices (\ref{eq-braiding}) and (\ref{eq-sigmatilde})
into the formulas $\sigma _{23}\ts _{12}\ts _{23}=
\ts _{12}\ts _{23}\sigma _{12}$ and $g_{23}\ts _{12}\ts _{23}=g_{12}$,
we obtain the requirements
\begin{gather*}
\epsfig{file=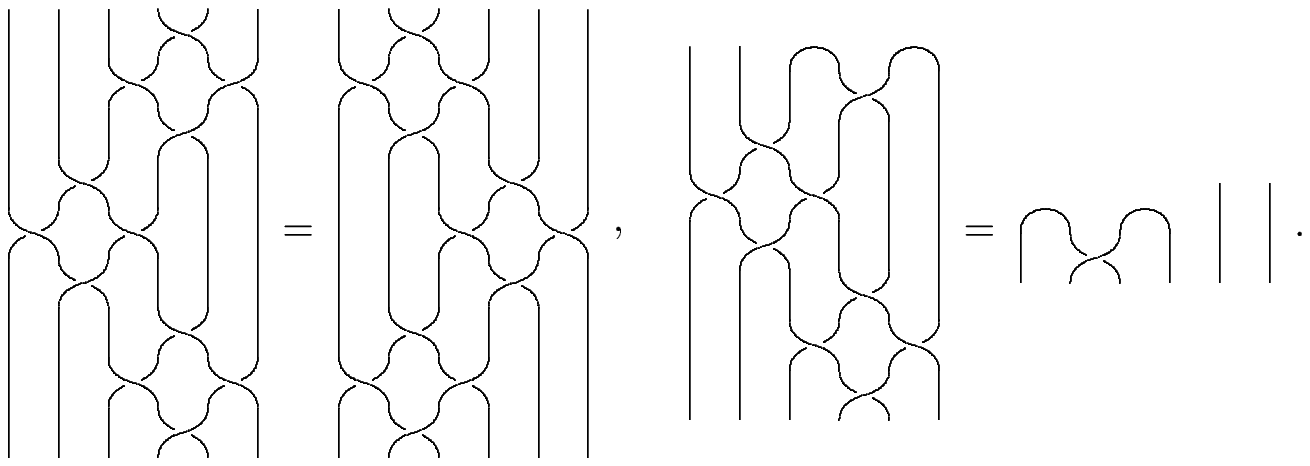}
\end{gather*}
But these formulas are also satisfied. Hence by Proposition \ref{s-clconn},
the linear left connection $\nabla $ on $\Gam $ given by $\nu \ts _{(1)}$,
$\nu \in \{1,-1\}$, 
is compatible with the multiplication of $\Cl (\Gam ,\sigma ,g)$.
\end{bew}

Now let us look for linear left connections $\szushg $ on the spinor module
$\mc{S}_0$ which are compatible with the right coaction of $\tilde{\A }$ on
$\mc{S}_0$.
This means that $\szushg $ should satisfy the equation
$\rkow (\szushg \psi )=(\szushg \ot \id )\rkow \psi $ for all
$\psi \in \mc{S}_0$.
Since $\Gam $ is an inner FODC, by Lemma \ref{l-innzushg} it
suffices to determine all homomorphisms
$\tau :\mc{S}_0\otA \Gam \to \Gam \otA \mc{S}_0$ and
$V:\mc{S}_0\to \Gam \otA \mc{S}_0$ of left-covariant bimodules,
which are compatible with $\rkow $.
Let $V:\mc{S}^\eta _0\to \Gam \otA \mc{S}^{\eta '}_0$,
$\eta ,\eta '\in \{+,-\}$, be a homomorphism of
left-covariant left $\A $-modules. Hence there exist $V^{ijk}_l\in \comp $,
$i,j,k,l=1,2$, such that
$V(\psi ^\eta _l)=V^{ijk}_l\w _{ij}\otA \psi ^{\eta '}_k$.
Compatibility of $V$ with $\rkow $ implies that
$(V^{ijk}_l)\in \Mor (u,u\ot u\ot u)$.
Using now (\ref{eq-psiuvert}) and (\ref{eq-wuvert}), one can determine
whether $V$ is a homomorphism of right $\A $-modules.
An equivalent criterion is that the matrix $\bar{V}:=(V^{ijk}_l)$
satisfies the equation
\begin{align}
\spinsig{0}\spinsig{\eta '}\Rda _{12}\Rdam _{23}\Rda ^{\eta '}_{34}\bar{V}_1&=
\spinsig{\eta }\bar{V}_2\Rda ^\eta _{12}\qquad \text{(leg numbering)}.
\end{align}
Since $\bar{V}\in \Mor(u,u\ot u\ot u)$, this implies that
\begin{align}\label{eq-HopfV}
V(\psi ^+_i)&=\alpha _+\oC ^{jk}\w _{ij}\otA \psi ^-_k,&
V(\psi ^-_i)&=\alpha _-\Rda {}^{kl}_{mi}\oC ^{jm}\w _{jk}\otA \psi ^+_l,
\end{align}
$\alpha _+,\alpha _-\in \comp $. Let $\bar{V}_\eta $, $\eta \in \{+,-\}$,
denote the complex matrix defined by
$V(\psi ^\eta _l)=(\bar{V}_\eta )_l^{ijk}\w _{ij}\otA \psi ^{-\eta }_k$.

Similar considerations lead to the determination of all homomorphisms $\tau $
of the left-covariant $\A $-bimodules $\mc{S}^\eta _0\otA \Gam $ and
$\Gam \otA \mc {S}^{\eta '}_0$, such that $\tau $ is compatible with $\rkow $.
Since $\tau $ is a homomorphism of left-covariant left $\A $-modules,
there exists a complex matrix $\bar{\tau }:=(\tau ^{rst}_{ijk})$,
such that $\tau (\psi ^\eta _i\otA \w _{jk})=q^{j+k}\tau ^{rst}_{ijk}
\w _{rs}\otA \psi ^{\eta '}_t$. Compatibility of $\tau $ with $\rkow $
implies that $\bar{\tau }\in \Mor (u\ot u\ot u)$. Moreover, $\tau $ is a
homomorphism of right $\A $-modules if and only if
\begin{align}
\spinsig{\eta }\bar{\tau }_{234}\Rda {}^\eta _{12}\Rda _{23}\Rdam _{34}&=
\spinsig{\eta '}\Rda _{12}\Rdam _{23}\Rda ^{\eta '}_{34}\bar{\tau }_{123}.
\end{align}
The solutions of this equation correspond to the mappings $\tau $ given by
\begin{equation}\label{eq-Hopftau}
\begin{aligned}
\tau (\psi ^+_i\otA \w _{jk})&=
q^{j+k}(\gamma ^+_1\delta ^r_i\Rda ^{st}_{jk}
+\gamma ^+_2\Rda ^{st}_{lk}\Rda ^{rl}_{ij})\w _{rs}\otA \psi ^+_t,\\
\tau (\psi ^-_i\otA \w _{jk})&=
q^{j+k}(\gamma ^-_1\Rdam {}^{rs}_{ij}\delta ^t_k
+\gamma ^-_2\Rda ^{st}_{lk}\Rdam {}^{rl}_{ij})\w _{rs}\otA \psi ^-_t,
\end{aligned}
\end{equation}
$\gamma ^\pm _1,\gamma ^\pm _2\in \comp $.
Let $\bar{\tau }_\eta $, $\eta \in \{+,-\}$,
denote the complex matrix defined by
$\tau (\psi ^\eta _i\otA \w _{jk})=
q^{j+k}(\bar{\tau }_\eta )_{ijk}^{rst}\w _{rs}\otA \psi ^\eta _t$.

Finally, to obtain a linear left connection on the spinor module $\mc{S}_0$
compatibility with the left action of $\Cl (\Gam ,\sigma ,g)$ has to hold.
This can be checked by Proposition \ref{s-modcompconn}. Observe that
because of equations (\ref{eq-HopfV}) and (\ref{eq-Hopftau}),
for $\psi \in \mc{S}_0^\eta $, $\eta \in \{+,-\}$, the first and second
summand of (\ref{eq-modcompconn}) is an element of
$\Gam \otA \mc{S}_0^{-\eta }$ and $\Gam \otA \mc{S}_0^\eta $, respectively.
Hence both summands have to vanish. Inserting formula (\ref{eq-sigmatilde})
for $\ts $ and equation (\ref{eq-clwirk})
and using the graphical calculus the following equations have to hold.

{\allowdisplaybreaks
\begin{gather}
\epsfig{file=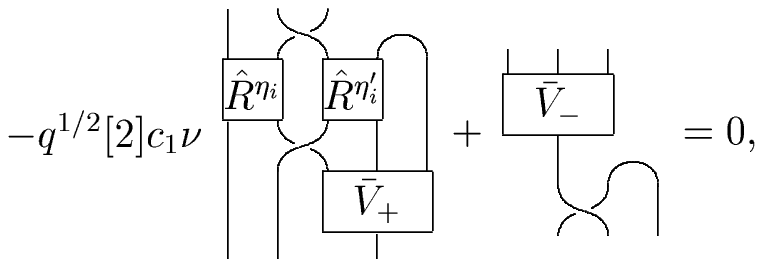}\\
\epsfig{file=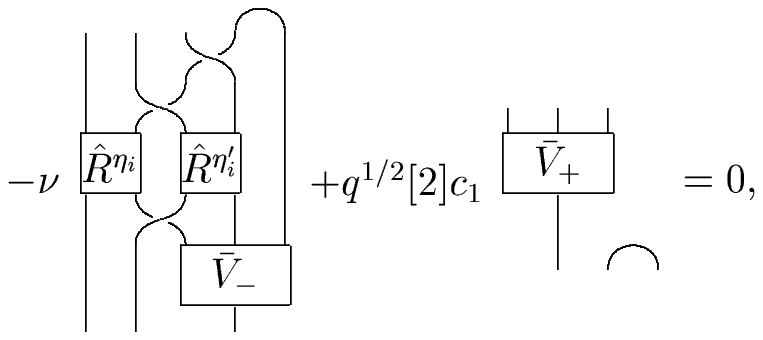}\\
\epsfig{file=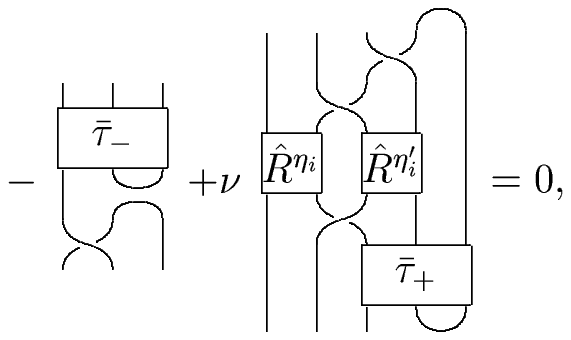}\\
\epsfig{file=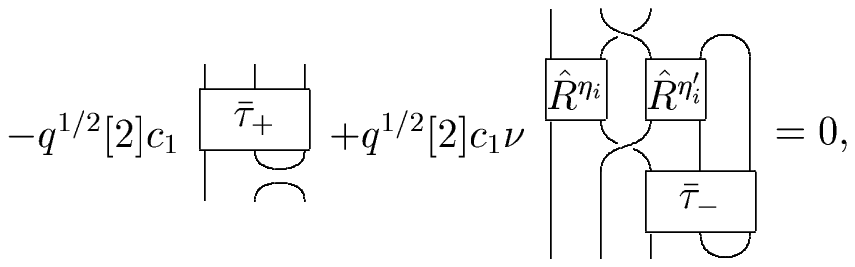}
\end{gather}
}
where $i\in \{1,2,3,4\}$ and $\nu \in \{+,-\}$.
We obtain the solutions $V=0$ and
\begin{equation}\label{eq-s0zushg}
\begin{aligned}
\tau (\psi _m^+\otA \w _{jk})&=\gamma q^{j+k}\Rda {}^{st}_{lk}
\Rda {}^{\eta _i}{}^{rl}_{mj}\w _{rs}\otA \psi ^+_t,\\
\tau (\psi _m^-\otA \w _{jk})&=\nu \gamma q^{j+k}\Rda {}^{\eta '_i}{}^{st}_{lk}
\Rdam {}^{rl}_{mj}\w _{rs}\otA \psi ^-_t,
\end{aligned}
\end{equation}
$\gamma \in \comp $. These considerations prove the first part of the
following theorem.

\begin{thm}\label{t-SLq2spinconn}
Let $g$ be one of the $\sigma $-metrics of Lemma \ref{l-bicsigmet} and let
$\nabla $ be one of the left connections in Proposition \ref{s-lconnslq2}.
Then there exists a 1-parameter family of linear left connections $\szushg $
on the spinor module $\mc{S}_0$. This connections are given by
(\ref{eq-innzushg}) with $V=0$ and $\ts :=\tau $ as in (\ref{eq-s0zushg}).
Moreover,\\
{\rm (i)} the left-covariant $\A $-subbimodules $\mc{S}_0^\eta $ of $\mc{S}_0$
are invariant with respect to $\szushg $ in the sense that
$\szushg \mc{S}_0^\eta \subseteq \Gam \otA \mc{S}_0^\eta $.\\
{\rm (ii)} the Dirac operator $D$ corresponding to $\szushg $ is symmetric
with respect to the metric $\hmet{\cdot }{\cdot }$ given in Theorem
\ref{t-unmet} if and only if
$q\in \real $ and $\bar{\gamma }=\gamma$, or $|q|=1$ and
$\bar{\gamma }=q^{6-3/2\eta _i-3/2\eta '_i}\nu \gamma $.
\end{thm}

For some computations it is necessary to know $\szushg $ in the more explicit form
\begin{equation}\label{eq-szushgSLq2}
\begin{aligned}
\szushg \psi ^+_m &=\left( \frac{q^{1/2}}{\qm }(1-q^{7/2-\eta _i/2}\gamma )\oC ^{jk}
\delta ^l_m-q^{5/2}\frac{1-\eta _i}{2}\gamma \delta ^j_m \oC ^{kl}\right)
\w _{jk}\otA \psi ^+_l,\\
\szushg \psi ^-_m &=\left( \frac{q^{1/2}}{\qm }(1-q^{7/2+\eta '_i/2}\nu \gamma )
\oC ^{jk}\delta ^l_m-q^{5/2}\frac{1+\eta '_i}{2}\nu \gamma \delta ^j_m \oC ^{kl}\right)
\w _{jk}\otA \psi ^-_l.
\end{aligned}
\end{equation}

\begin{bew}[ of (ii)]
Because of Proposition \ref{s-symmdirac} we only have to show that
$\hmet{D\varphi }{\psi }=\hmet{\varphi }{D\psi }$ for all
$\varphi ,\psi \in \linv{(\mc{S}_0)}$. {}From (\ref{eq-s0zushg}) and
(\ref{eq-clwirk}) we easily compute that
\begin{align}\label{eq-dpsi}
D\psi ^+_j&=\frac{q^{-1}}{\qm }(1-q^{9/2-3/2\eta _i}\gamma )\psi ^-_j,&
D\psi ^-_j&=\frac{-q\qp c_1}{\qm }(1-q^{3/2-3/2\eta '_i}\nu \gamma )\psi ^+_j.
\end{align}
Inserting this into (\ref{eq-unmet}) we obtain that $D$ is symmetric
if and only if $\gamma $ satisfies the given condition.
\end{bew}

\section{Invariant differential operators on $\mc{S}_0$}
\label{sec-invdo}

In Section \ref{sec-Dirac} it was shown that the Dirac operator and the
connection Laplacian on $\mc{S}_0$ are left-invariant differential operators.
In this section we will prove that they are compatible with the right
coaction $\rkow $ of $\tilde{\A }$ on $\mc{S}_0$, i.\,e.\
$\rkow (D\psi )=(D\ot \id )\rkow (\psi )$ and
$\rkow (\nabla ^*\szushg \psi )=(\nabla ^*\szushg \ot \id )\rkow (\psi )$
for $\psi \in \mc{S}_0$. Such left-invariant differential operators are called
\textit{invariant}.
On the other hand, sums, complex multiples and products of invariant
differential operators $\partial _1$ and $\partial _2$ on $\mc{S}_0$,
\begin{align}
(\partial _1+\partial _2)\psi &=\partial _1\psi +\partial _2\psi,&
(\lambda \partial _1)\psi &=\lambda (\partial _1\psi),&
(\partial _1\partial _2)\psi &=\partial _1(\partial _2\psi),
\end{align}
where $\lambda \in \comp $,
are again invariant differential operators on $\mc{S}_0$.
Hence the set of invariant differential operators on $\mc{S}_0$
forms an algebra $\mc{D}(\mc{S}_0)$.
Now we want to find out whether the Dirac operator and the connection
Laplacian are generic in the algebra $\mc{D}(\mc{S}_0)$ in some sense.

Let $\partial $ be a left-invariant differential operator on $\mc{S}_0$,
\begin{align}\label{eq-allglinvdo}
\partial (a\psi _i^\eta )&=\sum _{j,k,\eta '}a_{(1)}
p_{k,j}^{\eta ,\eta '}(a_{(2)})\uC _{ji}\psi _k^{\eta '}.
\end{align}
Then $\partial $ is invariant if and only if
$\rkow \partial (a\psi _i^\eta )=(\partial \ot \id )\rkow (a\psi _i^\eta )$
for $\eta \in \{+,-\}$, $i=1,2$, $a\in \A $.
This is equivalent to the equation
\begin{align}
a_{(1)}\ot a_{(2)}p_{k,j}^{\eta ,\eta '}(a_{(3)})\uC _{ji}u^l_k 
&=a_{(1)}\ot p_{l,k}^{\eta ,\eta '}(a_{(2)})\uC _{kj}a_{(3)}u^j_i,
\end{align}
$a\in \A $, $i,l=1,2$, $\eta ,\eta '\in \{+,-\}$.
Multiplying by $S(u^i_m)$ from the right, applying $S\ot \id $ and then
multiplying both factors of the tensor product we obtain
\begin{align}\label{eq-invdo}
p_{k,j}^{\eta ,\eta '}(a)\uC _{ji}u^l_k S(u^i_m)
&=S(a_{(1)})p_{l,j}^{\eta ,\eta '}(a_{(2)})\uC _{jm}a_{(3)},
\end{align}
$a\in \A $, $l,m=1,2$, $\eta ,\eta '\in \{+,-\}$.
In particular, $\partial $ is the sum of four invariant differential
operators $\partial _{\eta ,\eta '}:\mc{S}_0\to \mc{S}_0^{\eta '}$,
$\eta ,\eta '\in \{+,-\}$, such that
$\partial _{\eta ,\eta '}(\psi )=P_{\eta '}(\partial (\psi ))$
(projection to the component in $\mc{S}_0^{\eta '}$
with respect to the decomposition
$\mc{S}_0=\mc{S}_0^+\oplus \mc{S}_0^-$)
for $\psi \in \mc{S}_0^\eta $ and
$\partial _{\eta ,\eta '}(\psi )=0$ for $\psi \in \mc{S}_0^{-\eta }$.
{}From now on we suppose that $\partial =\partial _{+,+}$. The results carry
over to the other cases as well.

By Definition \ref{d-invdo} the functionals $p_{i,j}=p_{i,j}^{+,+}$
are elements of the unital subalgebra $\Xalg[\pm ]$ of $\Uqsl 2$ generated by
$\x _{kl}$, $k,l=1,2$. Evaluating a functional $f\in \Anull $ on both sides of
equation (\ref{eq-invdo}) yields that
\begin{align}
(\adR (f)p_{l,n})(a)&:=(S(f_{(1)})p_{l,j}f_{(2)})(a)\uC _{jm}\oC ^{mn}
\notag \\
\label{eq-adRfp}
&=p_{k,j}(a)\uC _{ji}f(u^l_k\oC ^{ir}\uC _{sm}u^s_r))\oC ^{mn}
=f(u^l_ku^n_r)p_{k,r}(a)
\end{align}
for all $l,n=1,2$, $f\in \Anull $, $a\in \A $.
Moreover, since $\Anull $ separates the elements of $\A $ (see
\cite[Theorem 11.22]{b-KS}), equation (\ref{eq-adRfp}) is equivalent to
(\ref{eq-invdo}).

Let $I(\partial )$ denote the pair $(v^0_0,(v^1_{-1},v^1_0,v^1_1))$, where
\begin{align}\label{eq-vmk}
v^m_k&:=\frac{-1}{\qp ^{1/2}}\sum _{i,j}C_q^{-1}(1/2,1/2,m;i-3/2,j-3/2,k)
p_{i,j},
\end{align}
$m=0,1$, $k=-m,\ldots ,m$,
and $C_q^{-1}(l_1,l_2,l;i,j,k)$ are the inverse Clebsch-Gordan coefficients
of the tensor product of the corepresentations $\coalg ^{(l_1)}$ and
$\coalg ^{(l_2)}$ (see Appendix \ref{sec-appslq2}). Then
\begin{align}\label{eq-uCGC}
\coalg ^{(l)}_{mn}C_q^{-1}(l_1,l_2,l;r,s,n)&=
\sum _{i,j}C_q^{-1}(l_1,l_2,l;i,j,m)\coalg ^{(l_1)}_{ir}\coalg ^{(l_2)}_{js}
\end{align}
and (\ref{eq-adRfp}) implies that
\begin{align}\label{eq-adRfv01}
\adR (f)v^0_0&=f(1)v^0_0,&
\adR (f)v^1_i&=f(\coalg ^{(1)}_{ij})v^1_j
\end{align}
for all $f\in \Anull $, $i=1,2,3$.

\begin{lemma}\label{l-idocorrv}
The mapping $I$ gives a one-to-one correspondence between invariant
differential operators $\partial $ on $\mc{S}_0^+$ and
pairs $v=(v^0_0,(v^1_{-1},v^1_0,v^1_1))$,
where $v^0_0,v^1_i\in \Xalg[\pm ]$ satisfy (\ref{eq-adRfv01}).
Under this correspondence we have
\begin{align}
I(\id )&=(\koun ,(0,0,0)),&
I(\lambda \partial +\lambda '\partial ')&=
\lambda I(\partial )+\lambda 'I(\partial ')
\end{align}
for $\lambda ,\lambda '\in \comp $ and
$\partial ,\partial '\in \mc{D}(\mc{S}_0^+)$. Moreover,
\begin{align}\label{eq-Ipartpart}
I(\partial \partial ')&=I(\partial '\partial )=
(v^0_0v'{}^0_0,(v^1_{-1}v'{}^0_0,v^1_0v'{}^0_0,v^1_1v'{}^0_0))
\end{align}
if $I(\partial )=(v^0_0,(v^1_{-1},v^1_0,v^1_1))$
and $I(\partial ')=(v'{}^0_0,(0,0,0))$.
\end{lemma}

\begin{bew}
Since $\coalg ^{(1/2)}\ot \coalg ^{(1/2)}\cong \coalg ^{(0)}\oplus
\coalg ^{(1)}$, $p_{i,j}$ can be reconstructed from $v$ with help of
(\ref{eq-vmk}) and the Clebsch-Gordan coefficients. Therefore $I$ is
bijective. Clearly, by equation (\ref{eq-allglinvdo}) the operator
$\partial =\id $ corresponds to $p_{i,j}=\oC ^{ij}\koun $.
{}From (\ref{eq-vmk}) and (\ref{eq-CGC}) we conclude that
$v^0_0=-1/\qp \uC _{ij}p_{i,j}=\koun $ and $v^1_i=0$ in the case
$\partial =\id $. Now we only have to prove (\ref{eq-Ipartpart}).
Since $v'{}^0_0=-1/\qp \uC _{ij}p'_{i,j}$ and $v'{}^1_k=0$, we get
$p'_{i,j}=\oC ^{ij}v'{}^0_0$. Hence $\partial '(a\psi ^+_k)=
a_{(1)}v'{}^0_0(a_{(2)})\psi ^+_k$ and
$\partial \partial '(a\psi ^+_k)=
a_{(1)}(p_{i,j}v'{}^0_0)(a_{(2)})\uC _{jk}\psi ^+_i$.
Equation (\ref{eq-adRfv01}) gives that $\partial '$ commutes with
$\partial $. Therefore (\ref{eq-Ipartpart}) follows from the definition of
$I$.
\end{bew}

By the general theory (see e.\,g.\ (14.39) in \cite{b-KS}) the equation
$\adR (f)\x _{ij}=f(u^i_ku^j_l)\x _{kl}$ holds. Since
$(\oC ^{ij})\in \Mor (1,u\ot u)$, it follows that
\begin{align}
\adR (f)(\x _{ij}+q^{1/2}/\qm \cdot \oC ^{ij})&=
f(u^i_ku^j_l)(\x _{kl}+q^{1/2}/\qm \cdot \oC ^{kl})
\end{align}
for $i,j=1,2$.
Let us define the functionals
$D^0_{ij}$ and $D^1_{ij}$, $i,j=1,2$, by
\begin{align}\label{eq-idogenpa}
D^0_{ij}&=-1/\qp \oC ^{ij}\uC _{kl}(\x _{kl}+q^{1/2}\oC ^{kl}/\qm ),&
D^1_{ij}&=\x _{ij}+q^{1/2}\oC ^{ij}/\qm -D^0_{ij}.
\end{align}
Then obviously $\adR (f)D^m_{ij}=f(u^i_ku^j_l)D^m_{kl}$ for $m=0,1$.
Hence the left-invariant differential operators $\partial _m$,
defined by
\begin{align}\label{eq-idogen}
\partial _m(a\psi _i^+)&:=a_{(1)}D^m_{kj}(a_{(2)})\uC _{ji}\psi _k^+,
\qquad m=0,1,
\end{align}
are invariant.

\begin{thm}\label{t-invdosp}
Recall the definitions (\ref{eq-idogenpa}) and (\ref{eq-idogen}).
The algebra of invariant differential operators $\mc{D}(\mc{S}_0^+)$ on
$\mc{S}_0^+$ is isomorphic to the free commutative unital algebra
$\tilde{\mc{D}}(\mc{S}_0^+)$ generated by
$\partial _0$ and $\partial _1$ and the relation
\begin{align}\label{eq-invdorel}
(\partial _1)^2+q^{-1}\qm \partial _0\partial _1-q^{-2}\partial _0^2=
-q^{-1}/\qm ^2\cdot \id .
\end{align}
\end{thm}

\begin{bew}
First we compute $\w ^2a\psi _i^+$ in two ways which will lead to
(\ref{eq-invdorel}). Of course, $\w ^2=g(\w ,\w )=-\qp c_1/\qm ^2$.
On the other hand, $\w a=\dif a+a\w =a_{(1)}(\x _{kl}+q^{1/2}\oC ^{kl}/\qm )
(a_{(2)})\w _{kl}$. Comparing $(\w (\w a))\psi _i^+$ and
$(\w \w )a\psi _i^+$, this and equations (\ref{eq-clwirk}) give
\begin{align}
&a_{(1)}\bigl((\x _{rs}+q^{1/2}\oC ^{rs}/\qm )(\x _{kl}+q^{1/2}\oC ^{kl}/\qm )
\bigr)(a_{(2)})\w _{rs}\w _{kl}\psi _i^+=\notag \\
&\quad =q^{1/2}\qp c_1\uC _{sm}\uC _{ni}\Rda {}^{mn}_{kl}a_{(1)}
\bigl((\x _{rs}+q^{1/2}\oC ^{rs}/\qm )(\x _{kl}+q^{1/2}\oC ^{kl}/\qm )\bigr)
(a_{(2)})\psi _r^+\notag \\
&\quad =-\qp c_1/\qm ^2\cdot a\psi _i^+.
\end{align}
Multiplying and setting $\Rda {}^{mn}_{kl}=q^{1/2}\delta ^m_k\delta ^n_l+
q^{-1/2}\oC ^{mn}\uC _{kl}$ the latter becomes
\begin{align*}\tag{$*$}
&a_{(1)}(\x _{rs}\uC _{sk}\x _{kl}\uC _{li}+q^{-1/2}\x _{rl}\uC _{li}\\
&\phantom{a_{(1)}(}+q^{-1}\x _{rl}\uC _{li}\uC _{mn}\x _{mn}
+q^{-1/2}/\qm \cdot \uC _{mn}\x _{mn}\delta ^r_i)(a_{(2)})\psi _r^+=0.
\phantom{)}
\end{align*}
Set
$\tilde{\partial }_0(a\psi _i^+):=a_{(1)}\uC _{rs}\x _{rs}(a_{(2)})\psi _i^+$
and
$\tilde{\partial }_1(a\psi _i^+):=a_{(1)}\x _{kj}(a_{(2)})\uC _{ji}\psi _k^+$.
Then
\begin{equation}\label{eq-parttrafo}
\begin{aligned}
\partial _0&=-1/\qp \tilde{\partial }_0+q^{1/2}/\qm \cdot \id ,\quad &
\partial _1&=\tilde{\partial }_1+1/\qp \tilde{\partial }_0,\\
\tilde{\partial }_0&=-\qp \partial _0+q^{1/2}\qp /\qm \cdot \id ,\quad &
\tilde{\partial }_1&=\partial _1+\partial _0-q^{1/2}/\qm \cdot \id .
\end{aligned}
\end{equation}
Moreover, ($*$) is equivalent to
\begin{align}
(\tilde{\partial }_1)^2+q^{-1/2}\tilde{\partial }_1
+q^{-1}\tilde{\partial }_1\tilde{\partial }_0
+q^{-1/2}/\qm \cdot \tilde{\partial }_0=0.
\end{align}
{}From this equation (\ref{eq-invdorel}) easily follows.

In \cite{a-HeckSchm2} was proved that
the locally finite part $\mc{F}(\Anull )$, defined by
\begin{align}
\mc{F}(\Anull )&:=\{f\in \Anull \,|\, \dim \adR (\Anull )f <\infty \},
\end{align}
is isomorphic to the vector space
\begin{align}
\mc{F}(\Anull )=\bigoplus _{n\in 1/2\mathbb{N}_0}\lfunc (\coalg ^{(n)})
\oplus \koun _-\bigoplus _{n\in 1/2\mathbb{N}_0}\lfunc (\coalg ^{(n)}).
\end{align}
Because of $\dif a=\w a-a\w $, $a\in \A $, we get $\x _{ij}=q^{1/2}/\qm
(\oC ^{kl}f^{kl}_{ij}-\oC ^{ij}\koun )$. Therefore
\begin{align*}
\kopr (\x _{ij}+q^{1/2}\oC ^{ij}/\qm \koun )&=
\x _{kl}\ot f^{kl}_{ij}
+\koun \ot \x _{ij}+q^{1/2}\oC ^{ij}/\qm \koun \ot \koun \\
&=(\x _{kl}+q^{1/2}\oC ^{kl}/\qm \koun )\ot f^{kl}_{ij}.
\end{align*}
This together with Theorem 4.1 in \cite{a-HeckSchm2} (or direct computation)
gives that $\koun _\pm \lfunc (\coalg ^{1/2})=
\Lin \{\x _{ij}+q^{1/2}\oC ^{ij}/\qm \koun \}$ for the $4D_\pm $-calculus.
Then Proposition 2.6 in \cite{a-HeckSchm2} gives that
\begin{align}
\Xalg[\pm ]=\bigoplus _{n\in 1/2\mathbb{N}_0}
\koun _\pm ^n\lfunc (\coalg ^{(n)}).
\end{align}
Moreover, $\adR (\Anull )(\koun _\pm \lfunc (\coalg ^{(n)}))\subset
\koun _\pm \lfunc (\coalg ^{(n)})$. Therefore $p_{i,j}$ can be written as
a finite sum of functionals $p^{(n)}_{i,j}$, $n\in 1/2\mathbb{N}_0$,
and each differential operator corresponding to the functionals
$p^{(n)}_{i,j}$ for a fixed $n\in 1/2\mathbb{N}_0$ is invariant.

For $n\in 1/2\mathbb{N}_0$ fix a nonzero matrix
$F(n)=(F(n)^i_j)$, $i,j=-n,-n+1,\ldots ,n$, such that
$F(n)\in \Mor (\coalg ^{(n)},(\coalg ^{(n)})\cont )$. Since $\coalg ^{(n)}$
is an irreducible corepresentation, $F(n)$ is invertible and unique up
to a complex factor. Now let us define the mappings
$\tlfunc :\coalg ^{(n)}\to \Anull $
by $\tlfunc (\coalg ^{(n)}_{ij})=\lfunc (\coalg ^{(n)}_{ik})F(n)^{-1}{}^j_k$.
Because of
$\adR (f)\lfunc (\coalg ^{(n)}_{ij})=
f(\coalg ^{(n)}_{ik}S(\coalg ^{(n)}_{lj}))\lfunc (\coalg ^{(n)}_{kl})$
we conclude that
\begin{align}
\adR (f)\tlfunc (\coalg ^{(n)}_{ij})&=
f(\coalg ^{(n)}_{ik}\coalg ^{(n)}_{jl})\tlfunc (\coalg ^{(n)}_{kl}).
\end{align}
{}From (\ref{eq-uCGC}) with $l_1=l_2=n$ we obtain the formula
and hence
\begin{align}
\adR (f)\tlfunc (C_q^{-1}(n,n,t;i,j,r)\coalg ^{(n)}_{ij})&=
f(C_q^{-1}(n,n,t;k,l,s)\coalg ^{(t)}_{rs})\tlfunc (\coalg ^{(n)}_{kl})
\end{align}
for all $r,n$ and all $t=0,1,\ldots ,2n$.
This proves that for each $n\in 1/2\mathbb{N}_0$ there exists a unique
1-dimensional complex subspace $V^{(n),1}$ of $\coalg ^{(n)}$
and a unique 3-dimensional one $V^{(n),3}$ with basis $\{v^{(n),1}_0\}$
and $\{v^{(n),3}_i\,|\,i=-1,0,1\}$, respectively, such that
\begin{align}\label{eq-adfv}
\adR (f)\tlfunc (v^{(n),2m+1}_i)&=
f(\coalg ^{(m)}_{ij})\tlfunc (v^{(n),2m+1}_j),& m=0,1,&\,
i=-m,\ldots ,m.
\end{align}
Moreover, this bases are unique up to a nonzero complex factor.

Setting $m=0$, (\ref{eq-adfv}) implies that $\tlfunc (v^{(n),1}_0)$ is a
central element in $\Anull $. In particular, $v^{(1/2),1}_0$ can be chosen
in such a way that
$\tlfunc (v^{(1/2),1}_0)$ becomes the quantum Casimir element $\Cas $
of $\Anull $. Recall that there is no nontrivial polynomial
function $p(\cdot )$ with complex coefficients such that
$p(\koun _\pm \Cas )\equiv 0$. Therefore,
since $\tlfunc (\coalg ^{(n_1)})\tlfunc (\coalg ^{(n_2)})=
\tlfunc (\coalg ^{(n_1)}\ot \coalg ^{(n_2)})$
for $n_1,n_2\in 1/2\mathbb{N}_0$,
$\koun _\pm ^{2n}\tlfunc (v^{(n),1}_0)$ is a polynomial function of
the element $\koun _\pm \Cas $ of degree $2n$. Further,
$\Cas \in \tlfunc (\coalg ^{(1/2)})$ implies that
$\Cas ^m\in \tlfunc ((\coalg ^{(1/2)})^{\ot m})$, $m\in \mathbb{N}$.
We obtain
\begin{gather}
\dim \bigoplus _{m=0}^n \koun ^{m+1}_\pm \Cas ^m \tlfunc (V^{(1/2),3})=3n+3,\\
\koun _\pm ^{m+1}\Cas ^m \tlfunc (V^{(1/2),3})\subset
\bigoplus _{{l=1\atop m+1-l\in 2\mathbb{Z}}}^{m+1} \koun ^{m+1}_\pm
\tlfunc (V^{(l/2),3})=
\bigoplus _{{l=1\atop m+1-l\in 2\mathbb{Z}}}^{m+1} \koun ^l_\pm
\tlfunc (V^{(l/2),3}),\\
\dim \bigoplus _{l=1}^{n+1} \koun ^l_\pm \tlfunc (V^{(l/2),3})=3n+3
\end{gather}
for all $n\in \mathbb{N}$. This means that for all $n\in \mathbb{N}$ there
exists a polynomial $P_n(\cdot )$ of degree $n-1$ such that
$\koun _\pm ^n \tlfunc (v^{(n/2),3}_i)=
P_n(\koun _\pm \Cas )\koun _\pm \tlfunc (v^{(1/2),3}_i)$, $i=-1,0,1$.
Denoting $I^{-1}\bigl((0,(
\koun _\pm ^n\tlfunc (v_{-1}^{(n/2),3}),
\koun _\pm ^n\tlfunc (v_0^{(n/2),3}),
\koun _\pm ^n\tlfunc (v_1^{(n/2),3})
))\bigr)$
by $\partial $, the latter equation and Lemma \ref{l-idocorrv} give
$\partial =P_n(\partial _0)\partial _1$.
Hence $\partial \in \tilde{\mc{D}}(\mc{S}_0^+)$.
\end{bew}

Let $(e_1,e_{-1})$ be a fixed basis of $\comp ^2$.

\begin{folg}\label{f-invdos0}
The algebra of invariant differential operators $\mc{D}(\mc{S}_0)$ on
$\mc{S}_0$ is isomorphic to the algebra $\tilde{\mc{D}}(\mc{S}_0):=
\comp ^2\ot \tilde{\mc{D}}(\mc{S}_0^+)\ot \comp ^2$
with multiplication
\begin{align}\label{eq-DS0mult}
(\alpha \ot \partial \ot \beta )(\alpha '\ot \partial '\ot \beta ')
&=\delta ^\beta _{\alpha '}\alpha \ot \partial \partial '\ot \beta ,
\end{align}
where $\alpha ,\alpha ',\beta ,\beta '\in \{e_1,e_{-1}\}$,
$\partial ,\partial '\in \tilde{\mc{D}}(\mc{S}_0^+)$.
\end{folg}

\begin{bew}
The mappings $\partial _{+-},\partial _{-+}:\mc{S}_0 \to \mc{S}_0$, defined by
$\partial _{+-}(\psi ^+_i)=\partial _{-+}(\psi ^-_i)=0$,
$\partial _{+-}(\psi ^-_i)=\psi ^+_i$, $\partial _{-+}(\psi ^+_i)=\psi ^-_i$,
are elements of $\mc{D}(\mc{S}_0)$.
One should identify the mapping
$e_1\ot \id \ot e_{-1}\in \tilde{\mc{D}}(\mc{S}_0)$
with the mapping $\partial _{+-}$ and
$e_{-1}\ot \id \ot e_1$ with the mapping $\partial _{-+}$, respectively.
Then the element $e_1\ot \partial '\ot e_1\in \tilde{\mc{D}}(\mc{S}_0)$ can be
identified with $\partial '\in \tilde{\mc{D}}(\mc{S}_0^+)$.
\end{bew}

By the defining equation (\ref{eq-dualzushg}) and by Theorem \ref{t-SLq2spinconn}
we are able to determine explicitly the left connection $\szushg ^*$ dual to $\szushg $.
Surprisingly the computations result in the formula
$\szushg ^*(\psi )=\w \otA \psi -\tau ^*(\psi \otA \w )$ for $\psi \in \mc{S}_0$, where
\begin{gather}\label{eq-s0dzushg}
\begin{aligned}
\tau ^*(\psi _m^+\otA \w _{jk})&=\nu \bar{\gamma }q^{j+k}\Rda {}^{st}_{lk}
\Rda {}^{-\eta '_i}{}^{rl}_{mj}\w _{rs}\otA \psi ^+_t,\\
\tau ^*(\psi _m^-\otA \w _{jk})&=\bar{\gamma }q^{j+k}\Rda {}^{-\eta _i}{}^{st}_{lk}
\Rdam {}^{rl}_{mj}\w _{rs}\otA \psi ^-_t,
\end{aligned}
\quad \text{for $q\in \real $,}\\
\begin{aligned}
\tau ^*(\psi _m^+\otA \w _{jk})&=q^{-6}\bar{\gamma }q^{j+k}\Rda {}^{st}_{lk}
\Rda {}^{\eta _i}{}^{rl}_{mj}\w _{rs}\otA \psi ^+_t,\\
\tau ^*(\psi _m^-\otA \w _{jk})&=\nu q^{-6}\bar{\gamma }q^{j+k}
\Rda {}^{\eta '_i}{}^{st}_{lk}\Rdam {}^{rl}_{mj}\w _{rs}\otA \psi ^-_t,
\end{aligned}
\quad \text{for $|q|=1$.}
\end{gather}
Comparing these formulas with (\ref{eq-s0zushg}) it turns out that
the left connection $\szushg ^*$ is a linear left connection on the spinor module
$\mc{S}_0$. However, for $q\in \real $, $\szushg $ and $\szushg ^*$ are compatible
with the same linear connection on $\Cl (\Gam ,\sigma ,g)$ if and only if
$\eta _i=-\eta '_i$, that is for $i=2$ and $i=3$.

\begin{folg}\label{f-Dinvop}
Let $\szushg $ be a linear connection on the spinor module $\mc{S}_0$.
Then the corresponding operators $D$ and $\nabla ^*\szushg $ are elements
of $\mc{D}(\mc{S}_0)$.
Moreover, both $D$ and $\nabla ^*\szushg $ are invariant first order
differential operators on $\mc{S}_0$.
\end{folg}

\begin{bew}
The assertion follows at once if we have shown that the formulas
\begin{align}\label{eq-diracslq2}
D&=e_{-1}\ot (-q^{1/2}\partial _1+q^{-3/2}\partial _0
-q^{7/2-3/2\eta _i}\gamma /\qm )\ot e_1\\
&\phantom{=}+e_1\ot (-q^{1/2}\qp c_1\partial _1-q^{1/2}\qp c_1\partial _0
+q^{5/2-3/2\eta '_i}\qp c_1\nu \gamma /\qm )\ot e_{-1},\notag \\
\label{eq-laplslq2}
\nabla ^*\szushg &=
e_1\ot (\alpha _1\partial _1+\alpha _0\partial _0+\alpha _2)\ot e_1\\
&\phantom{=}+e_{-1}\ot (\beta _1\partial _1+\beta _0\partial _0
+\beta _2)\ot e_{-1}\notag
\end{align}
hold, where
\begin{align*}
\alpha _1&=q^{7/2}c_1(\gamma (1-\eta _i)/2 +\nu \bar{\gamma }(1+\eta '_i)/2),\\
\alpha _0&=q^{5/2}c_1/\qm \cdot ((q^{3/2-\eta _i/2}+q^{-3/2+\eta _i/2})\gamma
+(q^{3/2+\eta '_i/2}+q^{-3/2-\eta '_i/2})\nu \bar{\gamma }),\\
\alpha _2&=-c_1/\qm ^2\cdot (\qp +q^6(q^{1+|\eta _i+\eta '_i|/2}
+q^{-1-|\eta _i+\eta '_i|/2})\nu \gamma \bar{\gamma }),\\
\beta _1&=q^{7/2}c_1(\bar{\gamma }(1-\eta _i)/2+\nu \gamma (1+\eta '_i)/2),\\
\beta _0&=q^{5/2}c_1/\qm \cdot ((q^{3/2-\eta _i/2}+q^{-3/2+\eta _i/2})\bar{\gamma }
+(q^{3/2+\eta '_i/2}+q^{-3/2-\eta '_i/2})\nu \gamma ),\\
\beta _2&=-c_1/\qm ^2\cdot (\qp +q^6(q^{1+|\eta _i+\eta '_i|/2}
+q^{-1-|\eta _i+\eta '_i|/2})\nu \gamma \bar{\gamma })\\
\intertext{for $q\in \real $ and}
\alpha _1&=q^{1/2}c_1(1-\eta _i)/2\cdot (q^3\gamma +q^{-3}\bar{\gamma }),\\
\alpha _0&=c_1/\qm \cdot (q^{5/2}(q^{3/2-\eta _i/2}+q^{-3/2+\eta _i/2})\gamma
+q^{-7/2}(q^{3/2-\eta _i/2}+q^{-3/2+\eta _i/2})\bar{\gamma }),\\
\alpha _2&=-c_1\qp /\qm ^2\cdot (1+\gamma \bar{\gamma }),\\
\beta _1&=q^{1/2}c_1(1+\eta '_i)/2\cdot (q^3\nu \gamma +q^{-3}\nu \bar{\gamma }),\\
\beta _0&=c_1/\qm \cdot (q^{5/2}(q^{3/2+\eta '_i/2}+q^{-3/2-\eta '_i/2})\nu \gamma
+q^{-7/2}(q^{3/2+\eta '_i/2}+q^{-3/2-\eta '_i/2})\nu \bar{\gamma }),\\
\beta _2&=-c_1\qp /\qm ^2\cdot (1+\gamma \bar{\gamma })
\end{align*}
for $|q|=1$.

By equation (\ref{eq-connection}) we obtain
$D(a\psi )=a_{(1)}X_{kl}(a_{(2)})\w _{kl}\psi +aD\psi $ for $a\in \A $,
$\psi \in \mc{S}_0$. Inserting (\ref{eq-dpsi}) and (\ref{eq-clwirk}) and
using the mappings $\tilde{\partial }_0$ and $\tilde{\partial }_1$ defined before
(\ref{eq-parttrafo}) this gives
\begin{align*}
D(a\psi _m^+)&=a_{(1)}\x _{kl}(a_{(2)})(-q^{1/2}\uC _{lm}\psi _k^-
-q^{-1/2}\uC _{kl}\psi _m^-)
+\frac{q^{-1}-q^{7/2-3/2\eta _i}\gamma }{\qm }a\psi _m^-\\
&=\left(e_{-1}\ot \left(-q^{1/2}\tilde{\partial }_1-q^{-1/2}\tilde{\partial }_0
+\frac{q^{-1}-q^{7/2-3/2\eta _i}\gamma }{\qm }\right)\ot e_1\right)a\psi _m^+.
\end{align*}
Similar computations for $D(a\psi _m^-)$ and the transformation formulas
(\ref{eq-parttrafo}) yield (\ref{eq-diracslq2}).

For the computation of (\ref{eq-laplslq2}) we use (\ref{eq-connlaplexp}). First,
the equation $g(\w \otA \w a)-g(\w \otA \w )a=0$ gives
\begin{align*}
a_{(1)}(\x _{kl}\x _{mn})(a_{(2)})g(\w _{kl},\w _{mn})
+2a_{(1)}\x _{kl}(a_{(2)})g(\w _{kl},\w )=0,
\end{align*}
and hence $a_{(1)}(\x _{kl}\x _{mn})(a_{(2)})g(\w _{kl},\w _{mn})=
-2q^{-1/2}c_1/\qm a_{(1)}\uC _{kl}\x _{kl}(a_{(2)})$.
Consider the case when $q\in \real $.
{}From (\ref{eq-szushgSLq2}) we derive that
\begin{align*}
&a_{(1)}\x _{kl}(a_{(2)})(g\ot \id )(\w _{kl}\otA (\szushg +\szushg ^*)(\psi ^+_m))=\\
&\qquad \qquad =\bigl(c_1/\qm \cdot (2q^{-1/2}-q^{2+\eta _i/2}\gamma 
-q^{2-\eta '_i/2}\nu \bar{\gamma })\tilde{\partial }_0\\
&\qquad \qquad \phantom{=\bigl(}
+q^{7/2}c_1((1-\eta _i)\gamma /2 +(1+\eta '_i)\nu \bar{\gamma }/2) 
\tilde{\partial }_1\bigr)(a\psi ^+_m) 
\end{align*}
and
\begin{align*}
&(g\ot \id )(\id \ot \szushg ^*)\szushg (\psi ^+_m)=
-c_1/\qm ^2\cdot \bigl(\qp -q^3(q^{3/2-\eta _i/2}+q^{-3/2+\eta _i/2})\gamma \\
&\quad -q^3(q^{3/2+\eta '_i/2}+q^{-3/2-\eta '_i/2})\nu \bar{\gamma }
+q^6(q^{1+|\eta _i+\eta '_i|/2}+q^{-1-|\eta _i+\eta '_i|/2})\nu \gamma \bar{\gamma }\bigr).
\end{align*}
These together with (\ref{eq-parttrafo}) give the first summand for the
connection Laplacian. The other one, corresponding to $\nabla ^*\szushg
{\upharpoonright }\mc{S}_0^-$, and the formula for $|q|=1$ can be determined
similarly.
\end{bew}

Let $q\in \real $. Suppose that $\szushg $ is one of the linear left
connections on $\mc{S}_0$ from Theorem \ref{t-SLq2spinconn} such that
$\nu =1$, $\bar{\gamma }=\gamma $ and $\eta _i=\eta '_i$ (that is $i=1$
or $i=4$). Then $D$ is symmetric by Theorem \ref{t-SLq2spinconn}.(ii), but
$\szushg $ and $\szushg ^*$ are not compatible with the same left connection
on $\Cl (\Gam ,\sigma ,g)$ (see the considerations before Corollary
\ref{f-Dinvop}).
For these connections one can formulate a modification of Bochner's theorem.

\begin{thm}
\label{t-SLq2Bochner}
Let $\szushg $ be as above and let $D$ and $\nabla ^*\szushg $
be the corresponding Dirac operator and connection Laplacian, respectively.
Then the operator
\begin{align}\label{eq-D2lapl}
(q+1)D^2-q^{1/2-3/2\eta _i}\qp \nabla ^*\szushg &=\eta _iq^{-\eta _i-1}\qp c_1
\frac{q^3-1}{\qm }\frac{q^6\gamma ^2-1}{\qm }\id
\end{align}
is an invariant differential operator of order zero.
\end{thm}

\begin{bew}
Using equations (\ref{eq-diracslq2}) and (\ref{eq-invdorel}) we obtain that
$D^2=e_1\ot x\ot e_1+e_{-1}\ot x\ot e_{-1}$ (in the notation of Corollary
\ref{f-invdos0}), where $x$ denotes the operator
\begin{equation}\label{eq-D2SLq2}
\begin{aligned}
\qp c_1/\qm \bigl( (q^{4-3/2\eta _i}-q^{3-3/2\eta '_i}\nu )\gamma \partial _1
&+(q^{4-3/2\eta _i}+q^{1-3/2\eta '_i}\nu )\gamma \partial _0\\
&-(1+q^{6-3/2\eta _i-3/2\eta '_i}\nu \gamma ^2)/\qm \bigr).
\end{aligned}
\end{equation}
Setting $\nu =1$ and $\eta '_i=\eta _i$, $x$ becomes
\begin{align*}
\qp c_1/\qm \bigl( q^{3-3/2\eta _i}(q-1)\gamma \partial _1
+q^{1-3/2\eta _i}(q^3+1)\gamma \partial _0
-(1+q^{6-3\eta _i}\gamma ^2)/\qm \bigr).
\end{align*}
On the other hand, inserting $\nu =1$, $\bar{\gamma }=\gamma $ and
$\eta '_i=\eta _i$ into equation (\ref{eq-laplslq2}) we directly obtain that
$\nabla ^*\szushg=e_1\ot y\ot e_1+e_{-1}\ot y\ot e_{-1}$, where $y$ denotes
the operator
\begin{align*}
q^{7/2}c_1\gamma \partial _1
+q^{1/2}c_1(q+1)(q^3+1)\gamma /\qm \partial _0
-c_1(\qp +q^6(q^2+q^{-2})\gamma ^2)/\qm ^2.
\end{align*}
{}From these formulas (\ref{eq-D2lapl}) immediately follows.
\end{bew}

Let us compute the eigenvalues of the Dirac operator $D$ corresponding to a
linear left connection on the spinor module $\mc{S}_0$. Since $D$ maps
$\mc{S}_0^\eta $ onto $\mc{S}_0^{-\eta }$, $\eta \in \{+,-\}$, it suffices to
determine the eigenvalues of $D^2$ on $\mc{S}_0^+$. Indeed, if
$\psi =\psi '\oplus \psi ''$, $\psi '\in \mc{S}_0^+$, $\psi ''\in \mc{S}_0^-$,
and $D\psi =\lambda \psi $, $\lambda \in \comp $, then
$D\psi '=\lambda \psi ''$ and $D\psi ''=\lambda \psi '$. Hence
$D^2\psi '=\lambda ^2 \psi '$. Conversely, if $D^2\psi '=\lambda ^2\psi '$,
$\psi '\in \mc{S}_0^+$, $\lambda \in \compx $, then
$\psi '\pm 1/\lambda ^{1/2}D\psi '$ is an eigenvector of $D$ to the eigenvalue
$\pm \lambda ^{1/2}$.

In all cases, the proof of the above theorem shows that
$D^2{\upharpoonright }\mc{S}_0^+$ is of the form
$\alpha _1\tilde{\partial }_1+\alpha _0\tilde{\partial }_0+\alpha _2\id $,
where $\alpha _1,\alpha _0,\alpha _2\in \comp $. Since $\mc{D}(\mc{S}_0^+)$ is
a commutative algebra, $\tilde{\partial }_0$ and $\tilde{\partial }_1$ have
common eigenvectors. Our first result will be that the linear hull of these
common eigenvectors span $\mc{S}_0^+$. Indeed, by left-invariance of
$\tilde{\partial }_0$ and $\tilde{\partial }_1$ the vector spaces
$V_n :=\coalg ^{(n)}\linv{(\mc{S}_0^+)}$, $n\in \mathbb{N}/2$, are invariant
under the action of these differential operators. On the other hand, 
$V_n$ splits under the right coaction $\rkow $ into the direct sum
$V_n=V_n^+\oplus V_n^-$, where
\begin{align*}
V_n^+&=\Lin \{\coalg ^{(n)}_rs\psi ^+_iC_q(n,1,n+1/2,s,i,t)\}\quad \text{and}\\
V_n^-&=\Lin \{\coalg ^{(n)}_rs\psi ^+_iC_q(n,1,n-1/2,s,i,t)\}\quad (V_0^-=\{0\})
\end{align*}
are non-isomorphic irreducible bicovariant vector spaces.
Since $\tilde{\partial }_0$ and $\tilde{\partial }_1$ are invariant differential
operators, by Schurs lemma they act by multiplication with a scalar on these
vector spaces. Because of (\ref{eq-slq2coreps}) we have
$\mc{S}_0^+=\bigoplus _{n\in \mathbb{N}/2}V_n$, and hence
$\tilde{\partial }_0$ and $\tilde{\partial }_1$ are commonly diagonalizable.

For the computation of the eigenvalues of
$\tilde{\partial }_0$ and $\tilde{\partial }_1$ on the eigenspaces $V^\pm _m$,
$m\in 1/2\mathbb{N}_0$,
we use the canonical embedding of $\coalg ^{(m)}$ into the $2m$-fold tensor
product of the corepresentation $\coalg ^{(1/2)}$. Let $P_{(2m)}$ denote the
unique covariant projection of $(\coalg ^{(1/2)})^{\ot 2m}$ onto its subcoalgebra
$\coalg ^{(m)}$. Of course we have
$\uC _{i,i+1}P_{(2m)}=P_{(2m)}\oC _{i,i+1}=0$ and
$\Rda _{i,i+1}P_{(2m)}=P_{(2m)}\Rda _{i,i+1}=q^{1/2}P_{(2m)}$
for $1\leq i<2m$.
Therefore the vector spaces $V^+_m$ and $V^-_m$ are generated by the elements
\begin{equation}\label{eq-Vgener}
\begin{gathered}
u^{j_1}_{l_1}u^{j_2}_{l_2}\cdots u^{j_{2m}}_{l_{2m}}\psi ^+_{l_{2m+1}}
P_{(2m+1)}{}^{l_1l_2\cdots l_{2m+1}}_{n_1n_2\cdots n_{2m+1}}
\quad \text{and}\\
u^{j_1}_{l_1}u^{j_2}_{l_2}\cdots u^{j_{2m}}_{l_{2m}}\psi ^+_{n_{2m+1}}
P_{(2m)}{}^{l_1l_2\cdots l_{2m}}_{n_1n_2\cdots n_{2m}}\oC ^{n_{2m}n_{2m+1}},
\end{gathered}
\end{equation}
respectively.

The definition $\dif a=\w a-a\w $, $\w=\sum _{i,j}q^{1/2}\oC ^{ij}/\qm \w _{ij}$,
and formula (\ref{eq-wuvert}) imply that
\begin{gather*}
\epsfig{file=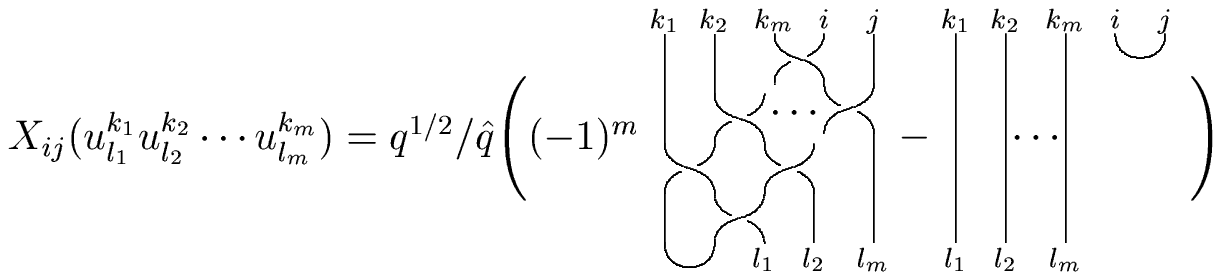}
\end{gather*}
Using the graphical calculus, from this one can easily compute the
eigenvalue of $\tilde{\partial }_0$ and $\tilde{\partial }_1$
to the eigenvectors given by (\ref{eq-Vgener}). We obtain that
\begin{align*}
\tilde{\partial }_0{\upharpoonright }V^+_m&=q^{1/2}/\qm \cdot
(q+q^{-1}-\nu _0^mq^{2m+1}-\nu _0^mq^{-2m-1})\id \quad (m\in 1/2\mathbb{N}_0),\\
\tilde{\partial }_0{\upharpoonright }V^-_m&=q^{1/2}/\qm \cdot
(q+q^{-1}-\nu _0^mq^{2m+1}-\nu _0^mq^{-2m-1})\id \quad (m\in 1/2\mathbb{N}),\\
\tilde{\partial }_1{\upharpoonright }V^+_m&=q^{1/2}/\qm \cdot
(\nu _0^mq^{2m}-1)\id \quad (m\in 1/2\mathbb{N}_0),\\
\tilde{\partial }_1{\upharpoonright }V^-_m&=q^{1/2}/\qm \cdot
(\nu _0^mq^{-2m-2}-1)\id \quad (m\in 1/2\mathbb{N}).
\end{align*}
{}From the transformation formulas (\ref{eq-parttrafo}) we get
\begin{align*}
\alpha _1\partial _1+\alpha _0\partial _0{\upharpoonright }V^+_m
&=\frac{q^{1/2}\nu _0^m}{\qp \qm }(\alpha _1(q^{2m-1}-q^{-2m-1})
+\alpha _0(q^{2m+1}+q^{-2m-1}))\id ,\\
\alpha _1\partial _1+\alpha _0\partial _0{\upharpoonright }V^-_m
&=\frac{q^{1/2}\nu _0^m}{\qp \qm }(\alpha _1(q^{-2m-3}-q^{2m+1})
+\alpha _0(q^{2m+1}+q^{-2m-1}))\id .
\end{align*}
Setting the correct parameter values for $D^2$ from (\ref{eq-D2SLq2})
we obtain the formulas
\begin{equation}
\begin{aligned}
D^2{\upharpoonright }V^+_m&=-\qp c_1\frac{q^{9/2+2m-3/2\eta _i}\nu _0^m\gamma
-1}{\qm }\frac{q^{3/2-2m-3/2\eta '_i}\nu _0^m\nu \gamma -1}{\qm }\id ,\\
D^2{\upharpoonright }V^-_m&=-\qp c_1\frac{q^{5/2-2m-3/2\eta _i}\nu _0^m\gamma
-1}{\qm }\frac{q^{7/2+2m-3/2\eta '_i}\nu _0^m\nu \gamma -1}{\qm }\id .
\end{aligned}
\end{equation}

\begin{appendix}
\section{On the quantum group $\SLq 2$}
\label{sec-appslq2}

Let $\OSLq 2$ denote the algebra generated by the elements $u^i_j$,
$i,j=1,2$, and relations $\Rda ^{ij}_{kl}u^k_mu^l_n=u^i_ku^j_l\Rda ^{kl}_{mn}$
and $u^1_1u^2_2-qu^1_2u^2_1=1$,
where $\Rda $ is given by (\ref{eq-RCmatrix}).
Setting $|u^i_j|=i-j$ it becomes a graded algebra with grading $|\cdot |$.
The coproduct $\kopr $ and the antipode $S$, where
\begin{align}
\kopr u^i_j&=u^i_k\ot u^k_j,&
S(u^i_j)&=\oC ^{ik}\uC _{lj}u^l_k
\end{align}
make $\OSLq 2$ to a Hopf algebra.
Further, $\OSLq 2$ is a coquasitriangular Hopf algebra with universal
$r$-form $\mathbf{r}$, such that $\mathbf{r}(u^i_j,u^k_l)=\Rda ^{ki}_{jl}$.

If $q$ is not a root of unity then the Hopf algebra $\OSLq 2$ is cosemisimple.
More precisely,
\begin{align}\label{eq-slq2coreps}
\OSLq 2=\bigoplus _{n\in \mathbb{N}_0/2}\coalg ^{(n)},
\end{align}
where $\coalg ^{(n)}$ is the linear span of the matrix elements of the unique
irreducible $2n+1$-dimensional corepresentation of $\OSLq 2$.
Let $\{\coalg ^{(n)}_{ij}\,|\,i,j=-n,-n+1,\ldots ,n\}$, $n\in \mathbb{N}_0/2$,
denote the basis of $\coalg ^{(n)}$ given in \cite{a-MMNNU91}.
Then we have $\kopr (\coalg ^{(n)}_{ij})=
\coalg ^{(n)}_{ik}\ot \coalg ^{(n)}_{kj}$.
Particularly,
\begin{align*}
(\coalg ^{(1/2)}_{ij})&=\begin{pmatrix}
u^1_1 & u^1_2\\
u^2_1 & u^2_2
\end{pmatrix},&
(\coalg ^{(1)}_{ij})&=\begin{pmatrix}
u^1_1u^1_1 & q'u^1_1u^1_2 & u^1_2u^1_2\\
q'u^1_1u^2_1 & 1+\qp u^1_2u^2_1 & q'u^1_2u^2_2\\
u^2_1u^2_1 & q'u^2_1u^2_2 & u^2_2u^2_2
\end{pmatrix},
\end{align*}
where $q'=(1+q^{-2})^{1/2}$.
Let $C_q(n_1,n_2,n;i,j,k)$ denote the Clebsch-Gordan coefficients of the
tensor product of the corepresentations $\coalg ^{(n_1)}$ and
$\coalg ^{(n_2)}$. Recall that
\begin{equation}\label{eq-CGC}
\begin{aligned}
\textstyle
C_q(\frac{1}{2},\frac{1}{2},0;i{-}\frac{3}{2},j{-}\frac{3}{2},0)&=
-1/\qp ^{1/2}\oC ^{ij}& \quad \text{for $i,j=1,2$,}\\
\textstyle
C_q^{-1}(\frac{1}{2},\frac{1}{2},0;i{-}\frac{3}{2},j{-}\frac{3}{2},0)&=
1/\qp ^{1/2}\uC _{ij}& \quad \text{for $i,j=1,2$.}
\end{aligned}
\end{equation}

The Hopf algebra $\OSLq 2$ carries three non-isomorphic $*$-structures
$\invcsu $, $\invncsu $ and $\invslR $ defined by
\begin{itemize}
\item
$q\in \mathbb{R}$;
$(u^1_1)^\invcsu =u^2_2$, $(u^1_2)^\invcsu =-qu^2_1$,
$(u^2_1)^\invcsu =-q^{-1}u^1_2$, $(u^2_2)^\invcsu =u^1_1$;
\item
$q\in \mathbb{R}$;
$(u^1_1)^\invncsu =u^2_2$, $(u^1_2)^\invncsu =qu^2_1$,
$(u^2_1)^\invncsu =q^{-1}u^1_2$, $(u^2_2)^\invncsu =u^1_1$;
\item
$|q|=1$;
$(u^i_j)^\invslR =u^i_j$ for $i,j=1,2$.
\end{itemize}
The Hopf $*$-algebras corresponding to the above involutions are denoted by
$\mc{O}(\mathrm{SU}_q(2))$, $\mc{O}(\mathrm{SU}_q(1,1))$ and
$\mc{O}(\mathrm{SL}_q(2,\mathbb{R}))$, respectively.

The functionals
$\lfunc (\coalg ^{(n)}_{ij})$ on $\OSLq 2$, defined by
\begin{align}
\lfunc (\coalg ^{(n)}_{ij})(a)&:=\mathbf{r}(\coalg ^{(n)}_{ik},a_{(1)})
\mathbf{r}(a_{(2)},\coalg ^{(n)}_{kj}),
\end{align}
are called the generalized $\ell $-functionals.

Let $\Uqsl 2$ denote the algebra generated by the elements $E,F,K$ and
$K^{-1}$ and relations
\begin{equation}
\begin{gathered}
KK^{-1}=K^{-1}K=1,\quad
KE=qEK,\quad KF=q^{-1} FK,\\
EF-FE=\frac{K^2-K^{-2}}{q-q^{-1}}.
\end{gathered}
\end{equation}
We denote the unit element in $\Uqsl 2$ by $\koun _+$.
In this paper we consider the central extension
$\tUqsl 2=\Uqsl2 \otimes \koun _-$ of $\Uqsl 2$.
The algebra $\tUqsl 2$ can be equipped with a Hopf structure such that
\begin{equation}
\begin{aligned}
\kopr (E)&=E\ot K+K^{-1}\ot E,\quad &\koun (E)&=0,\quad & S(E)&=-qE,\\
\kopr (F)&=F\ot K+K^{-1}\ot F,\quad &\koun (F)&=0, & S(F)&=-q^{-1}F,\\
\kopr (K)&=K\ot K, &\koun (K)&=1, & S(K)&=K^{-1},\\
\kopr (\koun _-)&=\koun _-\ot \koun _-, &\koun (\koun _-)&=1,&
S(\koun _-)&=\koun _-.
\end{aligned}
\end{equation}
There exists a dual pairing between the Hopf algebras $\tUqsl 2$ and $\OSLq 2$
such that for their generators the following formulas hold
(the matrix entries for $f\in \tUqsl 2$ are $f(u^i_j)$):
\begin{gather}\label{eq-pairing}
E=\begin{pmatrix}0&0\\1&0\end{pmatrix},\quad
F=\begin{pmatrix}0&1\\0&0\end{pmatrix},\quad
K=\begin{pmatrix}q^{-1/2}&0\\0&q^{1/2}\end{pmatrix},\quad
\koun _-=\begin{pmatrix}-1&0\\0&-1\end{pmatrix}.
\end{gather}
Hence $\tUqsl 2$ is a Hopf subalgebra of the dual Hopf algebra $\OSLq 2^\circ $
of $\OSLq 2$.

The real forms of $\OSLq 2$ induce $*$-structures on $\tUqsl 2$. They are
uniquely given by the pairing of $\OSLq 2$ and $\tUqsl 2$ and by
equation (\ref{eq-stern}). The explicit formulas are
\begin{itemize}
\item
$E^\invcsu =F$, $F^\invcsu =E$, $K^\invcsu =K$, $\koun _-^\invcsu =\koun _-$,
\item
$E^\invncsu =-F$, $F^\invncsu =-E$, $K^\invncsu =K$, $\koun _-^\invncsu =\koun _-$,
\item
$E^\invslR =E$, $F^\invslR =F$, $K^\invslR =K$, $\koun _-^\invslR =\koun _-$.
\end{itemize}

\section{Graphical calculus}
\label{sec-grcalc}

Here we collect the most important equations related to the matrices
$\Rda $, $\Rdam $, $\oC $ and $\uC $ in graphical form.
We use the symbols
\begin{gather*}
\begin{xy}
0;/r1.2pc/:
*{\mybox \Rda =}!R="mitte" +(.5,.5) \vcross ,
"mitte"+(1.6,0) *!L{\mybox ,}
\end{xy}
\quad
\begin{xy}
0;/r1.2pc/:
*{\mybox \Rdam =}!R="mitte" +(.5,.5) \vcrossneg ,
"mitte"+(1.6,0) *!L{\mybox ,}
\end{xy}
\quad
\begin{xy}
0;/r1.2pc/:
*{\mybox \oC =}!R="mitte" +(.5,.5) \vcap- ,
"mitte"+(1.6,0) *!L{\mybox ,}
\end{xy}
\quad
\begin{xy}
0;/r1.2pc/:
*{\mybox \uC =}!R="mitte" +(.5,-.5) \vcap ,
"mitte"+(1.6,0) *!L{\mybox,}
\end{xy}
\quad
\begin{xy}
0;/r1.2pc/:
*{\mybox \id =}!R="mitte"
+(.5,.5);p+(0,-1) **@{-} ,
"mitte"+(.6,0) *!L{\mybox .}
\end{xy}
\end{gather*}
{}From the settings (\ref{eq-RCmatrix}) it easily follows that
$\Rda {}^{ij}_{kl}=q^{1/2}\delta ^i_k\delta ^j_l+q^{-1/2}\oC ^{ij}\uC _{kl}$
and $\uC _{ij}\oC ^{ij}=-\qp $. This in turn implies the following formulas:
\begin{gather*}
\epsfig{file=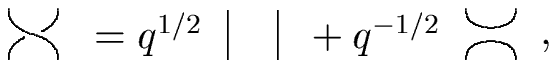}\quad
\epsfig{file=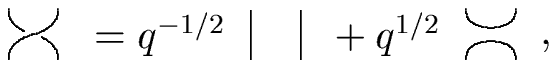}\\
\epsfig{file=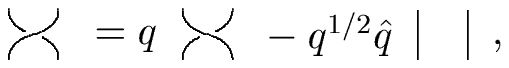}\quad
\epsfig{file=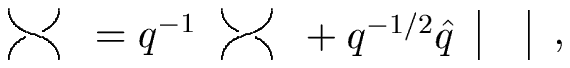}\\
\epsfig{file=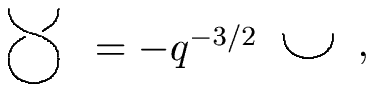}\quad
\epsfig{file=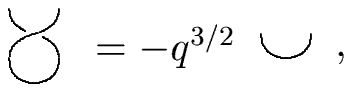}\quad
\epsfig{file=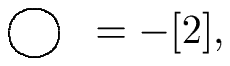}\\
\epsfig{file=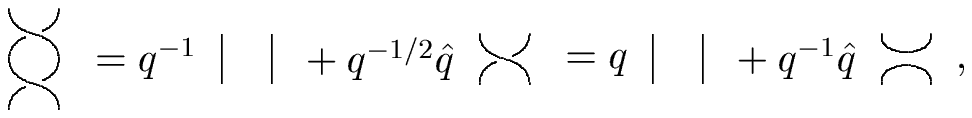}\\
\epsfig{file=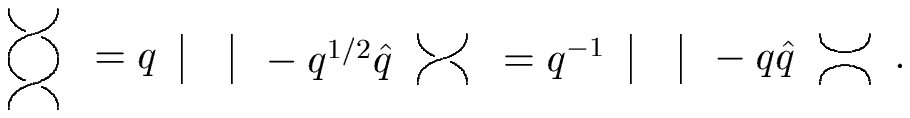}
\end{gather*}

\end{appendix}

\bibliographystyle{mybib}
\bibliography{quantum}

\end{document}